\pdfoutput=1
\RequirePackage{ifpdf}
\ifpdf 
\documentclass[pdftex]{sigma}
\else
\documentclass{sigma}
\fi

\usepackage{mathtools,tikz,tikz-cd,relsize,multirow}
\usepackage[all,cmtip]{xy}
\usepackage{blkarray}
\usetikzlibrary{positioning}
\usetikzlibrary{calc}
\usetikzlibrary{decorations.markings}
\usetikzlibrary{arrows}
\usetikzlibrary{calc}
\usetikzlibrary{decorations.markings}
\tikzstyle{hvector}=[inner sep=2pt,draw=blue!50,fill=blue!10,thick]
\tikzstyle{unit}=[inner sep=2pt,shape=circle, draw]
\tikzstyle{counit}=[inner sep=2pt,shape=circle, draw,fill=gray]
\tikzstyle{antipode}=[inner sep=2pt,shape=rectangle, draw]
\tikzstyle{cocycle}=[inner sep=2pt,shape=circle, draw]
\tikzstyle{twistedm}=[inner sep=2pt,shape=circle, fill=gray]
\tikzstyle{autom}=[inner sep=2pt,shape=circle, draw]
\tikzstyle{coact}=[inner sep=2pt,shape=circle, fill=black]

\makeatletter
\newif\ifv@ht@ \newif\ifv@dp@
\renewcommand*{\vphantom} {\v@ht@true \v@dp@true \h@false\ph@nt}
\renewcommand*{\hphantom} {\v@ht@false\v@dp@false\h@true \ph@nt}
\newcommand* {\htphantom}{\v@ht@true \v@dp@false\h@false\ph@nt}
\newcommand* {\dpphantom}{\v@ht@false\v@dp@true \h@false\ph@nt}
\renewcommand*{\phantom} {\v@ht@true \v@dp@true \h@true \ph@nt}
\renewcommand*{\finph@nt}{%
 \setbox2=\null
 \ifv@ht@ \ht2=\ht0 \fi
 \ifv@dp@ \dp2=\dp0 \fi
 \ifh@ \wd2=\wd0 \fi
 \box2 %
}
\makeatother

\numberwithin{equation}{section}

\newtheorem{Theorem}{Theorem}[section]
\newtheorem*{Theorem*}{Theorem}

\newtheorem{Lemma}[Theorem]{Lemma}

\theoremstyle{definition}
\newtheorem{Definition}[Theorem]{Definition}

\newtheorem{Remark}[Theorem]{Remark}

\newcommand{\matindex}[1]{\mbox{\scriptsize#1}}

\makeatletter
\def\namedlabel#1#2{\begingroup
 #2%
 \def\@currentlabel{#2}%
 \phantomsection\label{#1}\endgroup
}
\makeatother

\begin{document}

\allowdisplaybreaks

\newcommand{\arXivNumber}{2505.01398}

\renewcommand{\PaperNumber}{050}

\FirstPageHeading

\ShortArticleName{Extending Knot Polynomials of Braided Hopf Algebras to Links}

\ArticleName{Extending Knot Polynomials \\ of Braided Hopf Algebras to Links}

\Author{Stavros GAROUFALIDIS~$^{\rm a}$, Matthew HARPER~$^{\rm b}$, Ben-Michael KOHLI~$^{\rm c}$, Jiebo SONG~$^{\rm d}$\newline and Guillaume TAHAR~$^{\rm d}$}

\AuthorNameForHeading{S.~Garoufalidis, M.~Harper, B.-M.~Kohli, J.~Song and G.~Tahar}

\Address{$^{\rm a)}$~International Center for Mathematics, Department of Mathematics,\\
\hphantom{$^{\rm a)}$}~Southern University of Science and Technology, Shenzhen, P.R.~China}
\EmailD{\mail{stavros@mpim-bonn.mpg.de}}
\URLaddressD{\url{http://people.mpim-bonn.mpg.de/stavros/}}

\Address{$^{\rm b)}$~Department of Mathematics, Michigan State University, East Lansing, MI~48824, USA}
\EmailD{\mail{mrhmath@proton.me}, \mail{harpe111@msu.edu}}
\URLaddressD{\url{https://mrhmath.github.io/}}

\Address{$^{\rm c)}$~Section de Math\'ematiques, Universit\'e de Gen\`eve, \\
\hphantom{$^{\rm c)}$}~rue du Conseil-G\'en\'eral 7-9, 1205 Gen\`eve, Switzerland}
\EmailD{\mail{bm.kohli@protonmail.ch}}

\Address{$^{\rm d)}$~Beijing Institute of Mathematical Sciences and Applications, Beijing, P.R.~China}
\EmailD{\mail{songjiebo@bimsa.cn}, \mail{guillaume.tahar@bimsa.cn}}

\ArticleDates{Received October 20, 2025, in final form May 10, 2026; Published online May 18, 2026}

\Abstract{Recently, a plethora of multivariable knot polynomials were introduced by Kashaev and one of the authors, by applying the Reshetikhin--Turaev functor to rigid $R$-matrices that come from braided Hopf algebras with automorphisms. We~study the extension of these knot invariants to links, and use this to identify some of them with known link invariants, as conjectured in that same recent work.}


\Classification{57K14; 57K16; 17B37; 18M15}

\section{Introduction}
\label{sec:intro}

\subsection{From knot polynomials to link polynomials}
Suppose one is given a polynomial invariant of knots in 3-space. Is there a
natural way to extend it to a polynomial invariant of links? Of course, one
can extend it by declaring it to vanish for links of more than one component,
or that it is the product of the polynomials of each component of the link.
Most people would accept, however, that such extensions are not very natural.

This paper aims to give an answer to this question for the polynomial invariants
of (long) knots constructed via the Reshetikhin--Turaev functor~\cite{RT}
using as input a rigid $R$-matrix, as explained in detail in~\cite{Kas1}.
A rich source of rigid $R$-matrices was recently discovered by Kashaev and one
of the authors in~\cite{GK:multi} to come from braided Hopf algebras (such as Nichols
algebras) with automorphisms and their finite-dimensional left/right
Drinfeld--Yetter modules. This, together with the classification of Nichols algebras
gives a systematic way to construct multivariable polynomials of long knots.

This construction provides a novel framework to study quantum group invariants even if the invariants they produce are already known. For example, it provides a somewhat unified setting to consider the colored Jones polynomial and ADO polynomial invariants at roots of unity~\cite{ADO} from rank one Nichols algebras. One exchanges the structure of quantum $\mathfrak{sl}_2$ modules for the structure of a Nichols algebra, where the underlying vector spaces are naturally identified.
One also finds that the $\mathfrak{gl}(n|1)$ quantum invariants, i.e., the higher rank Links--Gould invariants, are associated to exterior algebras \cite{Kas3}.
Finally, the approach of \cite{GK:multi} implies a natural hierarchy of invariants. In the above reference,
2-variable knot polynomials coming from the simplest rank~2 Nichols algebra were
defined, and denoted by $\Lambda_\omega(t,s)$ and $V_n(t,q)$, respectively. Our goal is to extend these two 2-variable polynomial invariants of knots to invariants of links, especially towards applications in our sequel works \cite{GHKKST, GHKKW}.

Throughout the paper, all knots and links will be oriented, embedded in $S^3$,
and considered up to ambient isotopy.
But even more, a rigid $R$-matrix $R \in \mathrm{End}(V \otimes V)$ on a~vector space~$V$ together with
an enhancement $h \in \mathrm{End}(V)$ that satisfies the polynomial equations
from~\cite[Theorem~3.7]{Ohtsuki} gives invariants of tangles (up to isotopy)
that may have closed components. We review these identities in~\eqref{eq1} of Section~\ref{sub:tangles} and show in Lemma~\ref{lem.V10} that the rigid $R$-matrix for~$V_1$, introduced in~\cite{GK:multi}, has a canonically defined enhancement that satisfies the
required identities. The construction of this enhancement amounts to solving a system of equations, and a similar derivation applies to the $R$-matrix enhancements of $\Lambda_1$ and $\Lambda_{-1}$.

In this context, by canonical, we mean that upon assigning a natural (weight) grading to vectors and assuming that $h$ respects this grading (it is diagonal), then the relations \eqref{eq1} uniquely determine $h$. See Lemma~\ref{lem.V10}, for example, which uses \eqref{eq1b} to determine $h$. The sense in which $h$ is canonical can also be compared to the case of quantum groups equipped with an $R$-matrix. The ribbon element, if one exists, is unique up to a group-like element of order~2. This determines $h$ up to some signs which are fixed by the Reidemeister I constraints of \eqref{eq1b}.

Assuming this, a rigid $R$-matrix on $V$ gives rise to an invariant $F$ of
closed links. However, for some rigid $R$-matrices (e.g., for all discussed
in our paper) $F$ is the zero invariant.
To obtain a nontrivial invariant, we consider tangles,
where a rigid $R$-matrix gives an $\mathrm{End}(V)$-valued invariant $F_T$ of a
$(1,1)$-tangle $T$ and an $\mathrm{End}(V\otimes V)$-valued invariant $F_T$
of a $(2,2)$-tangle with upward oriented boundary, by which we mean a tangle
with two upward incoming and outgoing arcs and perhaps additional closed components.
To get to a scalar-valued invariant of links, we need to ensure:
\begin{enumerate}\itemsep=0pt
\item[(\namedlabel{item:P1}{$P_1$})]
 For every $(1,1)$-tangle $T$, $F_T$ is a scalar multiple of
 $\mathrm{id}_V$.
\item[(\namedlabel{item:P2}{$P_2$})]
 For every $(2,2)$-tangle $T$ defining a map $F_T\in\mathrm{End}\bigl(V^{\otimes 2}\bigr)$
 with left and right closures $T_1$ and $T_2$, we have
 $F_{T_1}=F_{T_2} \in \mathrm{End}(V)$.
\end{enumerate}

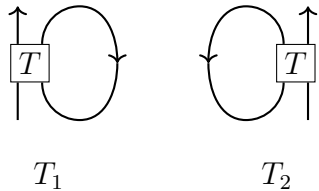
\begin{figure}[htpb!]\centering
\begin{tikzpicture}[baseline=10, xscale=1, yscale=1]
\node[rectangle,inner sep=3pt,draw] (c) at (1/2,1/2) {$T$};
\draw[thick,->] ($(c.north east)+(-.1,0)$) to [out=90, in= 180]
($(c.north east)+(.4,1/2)$) to [out=0, in= 90] ($(c.east)+(.9,0)$);
\draw[thick] ($(c.east)+(.9,0)$) to [out=-90, in= 0] ($(c.south east)+(.4,-1/2)$)
to [out=180, in= -90] ($(c.south east)+(-.1,0)$);
\draw[thick] ($(c.south west)+(.1,-1/2)$) to ($(c.south west)+(.1,0)$);
\draw[thick,->] ($(c.north west)+(.1,0)$) to ($(c.north west)+(.1,1/2)$);
\node[] at ($(c)+(.25,-1.5)$) {$T_1$};
\end{tikzpicture}
\qquad
\begin{tikzpicture}[baseline=10, xscale=1, yscale=1]
\node[rectangle,inner sep=3pt,draw] (c) at (1/2,1/2) {$T$};
\draw[thick,->] ($(c.north west)+(.1,0)$) to [out=90, in= 0]
($(c.north west)+(-.4,1/2)$) to [out=180, in= 90] ($(c.west)+(-.9,0)$);
\draw[thick] ($(c.west)+(-.9,0)$) to [out=-90, in= 180] ($(-.4,-1/2)+(c.south west)$)
to [out=0, in= -90] ($(c.south west)+(.1,0)$);
\draw[thick] ($(c.south east)+(-.1,-1/2)$) to ($(c.south east)+(-.1,0)$);
\draw[thick,->] ($(c.north east)+(-.1,0)$) to ($(c.north east)+(-.1,1/2)$);
\node[] at ($(c)+(-.25,-1.5)$) {$T_2$};
\end{tikzpicture}
\caption{The left and right closures $T_1$ and $T_2$ of a $(2,2)$-tangle
 $T$.}
\label{fig:T12}
\end{figure}

When \eqref{item:P1} is satisfied and $T$ is a $(1,1)$-tangle, we write
$F_T=\langle F_T\rangle \mathrm{id}_V$. Properties~\eqref{item:P1} and~\eqref{item:P2}
imply the existence of a scalar invariant of links $L$ obtained by cutting a link
along any component to obtain a $(1,1)$-tangle $L^\mathrm{cut}$.
We denote this modified Reshetikhin--Turaev link invariant
$\text{mRT}_L=\langle F_{L^\mathrm{cut}}\rangle$.

\begin{Remark}
Suppose that $F$ is determined by the monoidal category generated by the objects~$V$ and~$V^*$, and by the morphisms $R\in\mathrm{End}(V\otimes V)$ and the evaluation and coevaluation maps of Figure~\ref{fig:FIG3}
 defined by equations \eqref{eq1}. Then the properties~\eqref{item:P1} and~\eqref{item:P2} are the assumptions that the object $V$ is simple and that $V$ is an ambidextrous object in the sense of~\cite{GeerKujawaPatureau, GPT}. In this case, $\text{mRT}$ is known to define a link invariant by \cite[Theorem~3]{GPT}.
\end{Remark}

Once the link invariant is defined, the next question is how to identify it
with a known link invariant. For this, observe that any identity of the $R$-matrix
of braids implies a corresponding skein theory for the constructed link invariants.

A key example of this exposition is the $R$-matrix of the Alexander polynomial
given in the beautiful exposition of Ohtsuki~\cite[Proposition~2.6]{Ohtsuki}, which we
review in Section~\ref{sub:alex} below, using a~slightly different $R$-matrix for
the purpose of our work.

\subsection{Our results}
\label{sub:results}

With the above preliminaries, we can phrase our results.

\begin{Theorem}\label{thm.1}
The knot polynomial invariants $V_1$, $\Lambda_1$ and $\Lambda_{-1}$ extend
to invariants of links.
\end{Theorem}

The next theorem identifies the associated link invariants
$\Lambda_1$ and $\Lambda_{-1}$ confirming a conjecture in~\cite{GK:multi}.

\begin{Theorem}\label{thm.2}
For all links $L$, we have
\begin{gather}
\Lambda_{1,L}(t_0, t_1) = \Delta_L (t_0) \Delta_L(t_1)
\in \mathbb{Z}\bigl[t_0^{\pm 1},t_1^{\pm 1}\bigr],
\nonumber\\
\label{Lam1}
\Lambda_{-1,L}\bigl(t^{-2}, s^{-2}\bigr) = \Delta_{\mathfrak{sl}_3,L}(t,s)
\in \mathbb{Z}\bigl[t^{\pm2},s^{\pm2}\bigr],
\end{gather}
where $\Delta$ is the $($palindromic normalization of the$)$ Alexander polynomial and
$\Delta_{\mathfrak{sl}_3}$ is the invariant studied in~{\rm \cite{Harper}}.
\end{Theorem}

The Garoufalidis--Kashaev construction seems to offer a more natural setting for
connecting quantum invariants to invariants from classical topology. Specifically,
the coefficients of the $R$-matrices are valued in integral Laurent polynomials with
canonical variables; $t_1$, $t_2$ rather than their squares.

In subsequent work \cite{GHKKST}, we prove that the extension of $V_1$ coincides with the
Links--Gould invariant of links, verifying another conjecture of~\cite{GK:multi}. Building on both of these results, we prove in \cite{GHKKW} that the Garoufalidis--Kashaev invariant $V_2$ agrees with the 2-colored Links--Gould polynomial and we use this result to prove a conjecturally sharp bound for the 3-genus of a knot.

We should mention that the particular choice of local, tangle definition of these multivariable polynomials
of links plays an important role in their efficient computation, as was explained
time and again by Bar-Natan and van der Veen~\cite{BNV:PG,BNV:API}.
In particular, it leads to an efficient computation of the $V_n$ and the
$\Lambda_\omega$ polynomials of knots, that itself leads to newly discovered
patterns in knot theory; see~\cite{GL:patterns}.

\subsection{Plan of the proofs}

In Section~\ref{sec:RT}, we recall the basics on tangles and the Reshetikhin--Turaev
functor, as well as the properties of enhanced $R$-matrices that are
used to realize the functor, and its comparison with the tangle construction of~\cite{GK:multi}. We also introduce the notion of a weak conjugacy between pairs of enhanced $R$-matrices. Following Ohtsuki, we discuss in detail the
Alexander $R$-matrix of a 2-dimensional space and how it leads to a link invariant,
which is then identified with the Alexander polynomial.

We then discuss the rigid $R$-matrix of $V_1$ on a 4-dimensional space in Section~\ref{sec:V1} and produce an enhancement. We show that the enhanced $R$-matrix
defines a link invariant.

The extension of $\Lambda_1$ to a link invariant and its identification with a product
of Alexander polynomials is done simultaneously. This is proven by an appropriate
conjugation of the enhanced $R$-matrix at the level of individual tensor factors $V$
on which the $R$-matrix acts.

The proof of Theorems~\ref{thm.1} and~\ref{thm.2} for $\Lambda_{-1}$ is also based on a conjugacy of enhanced $R$-matrices, but is instead defined on the vector spaces $V^{\otimes n}$ for each $n\geq 1$. This determines an equivalence of braid group representations that respects partial trace operations, and implies equality with the link invariant $\Delta_{\mathfrak{sl}_3}$.

\section{Link invariants from the Reshetikhin--Turaev functor}
\label{sec:RT}

\subsection{A review of the Reshetikhin--Turaev functor}
\label{sub:RT}

In this section we review briefly the well-known Reshetikhin--Turaev functor
~\cite{RT,Tu:book} from tangles to tensor products of endomorphisms of a vector
space and its dual. All tangles are oriented and framed, however our invariants
will be independent of the framing.

\begin{figure}[htpb!]\centering
\begin{tikzpicture}[xscale=0.4, yscale =0.4]
\draw[thick,->] (0,0) to (0,2);
\draw[thick,->] (2,2) to (2,0);
\draw[thick,->] (4,0) to [out=90, in=180] (5,1);
\draw[thick,] (5,1) to [out=0, in=90] (6,0);
\draw[thick,->] (10,0) to [out=90, in=0] (9,1);
\draw[thick,] (9,1) to [out=180, in=90] (8,0);
\draw[thick,->] (14,2) to [out=-90, in=0] (13,1);
\draw[thick,] (13,1) to [out=180, in=-90] (12,2);
\draw[thick,->] (16,2) to [out=-90, in=180] (17,1);
\draw[thick,] (17,1) to [out=0, in=-90] (18,2);
\draw[thick,->] (22,0) to [out=90, in=-90] (20,2);
\draw[line width=3pt, white] (20,0) to [out=90, in=-90] (22,2);
\draw[thick,->] (20,0) to [out=90, in=-90] (22,2);
\draw[thick,->] (24,0) to [out=90, in=-90] (26,2);
\draw[line width=3pt, white] (26,0) to [out=90, in=-90] (24,2);
\draw[thick,->] (26,0) to [out=90, in=-90] (24,2);
\end{tikzpicture}
\caption{Elementary oriented tangles.}
\label{fig:tangles1}
\end{figure}
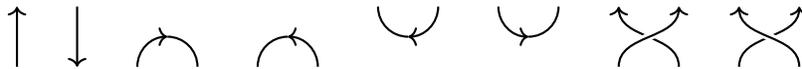

A tangle diagram can be decomposed into a finite number of elementary
oriented tangle diagrams shown in Figure~\ref{fig:tangles1}. Moreover, up to
 isotopy, we may assume that for any horizontal line drawn on the diagram there is
 at most one critical point or a single crossing.

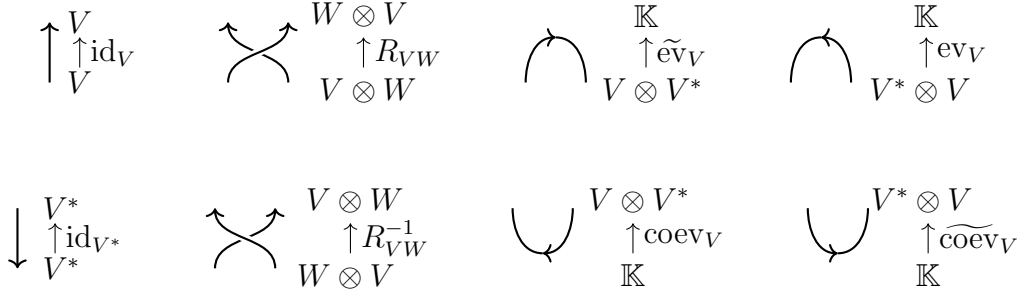
\begin{figure}[htpb!]\centering
\begin{tikzpicture}[baseline=0, xscale=.4, yscale=.4]
\draw[thick,->] (0,0) to (0,2);
\node at (1,0) {$V$};
\node at (1,2) {$V$};
\draw[->] (1,.5) to (1,1.5);
\node[right] at (1,1) {$\mathrm{id}_V$};
\end{tikzpicture}
\qquad
\begin{tikzpicture}[baseline=0, xscale=.4, yscale=.4]
\draw[thick,->] (7,0) to [out=90, in=-90] (5,2);
\draw[line width=3pt, white] (5,0) to [out=90, in=-90] (7,2);
\draw[thick,->] (5,0) to [out=90, in=-90] (7,2);
\node[below] at (9.5,.5) {$\phantom{W}V\otimes W\phantom{V}$};
\node[above] at (9.5,1.5) {$\phantom{V}W\otimes V\phantom{W}$};
\draw[->] (9.5,.5) to (9.5,1.5);
\node[right] at (9.5,1) {$R_{VW}$};
\end{tikzpicture}
\qquad
\begin{tikzpicture}[baseline=0, xscale=.4, yscale=.4]
\draw[thick,->] (4,0) to [out=90, in=180] (5,1.5);
\draw[thick,] (5,1.5) to [out=0, in=90] (6,0);
\node[below] at (8,.5) {$\phantom{^*}V\otimes V^*$};
\node[above] at (8,1.5) {$\mathbb{K}$};
\draw[->] (8,.5) to (8,1.5);
\node[right] at (8,1) {$\overrightarrow{\operatorname{ev}}_V$};
\end{tikzpicture}
\qquad
\begin{tikzpicture}[baseline=0, xscale=.4, yscale=.4]
\draw[thick,->] (10,0) to [out=90, in=0] (9,1.5);
\draw[thick,] (9,1.5) to [out=180, in=90] (8,0);
\node[below] at (12.5,.5) {$V^*\otimes V\phantom{^*}$};
\node[above] at (12.5,1.5) {$\mathbb{K}$};
\draw[->] (12.5,.5) to (12.5,1.5);
\node[right] at (12.5,1) {$\overleftarrow{\operatorname{ev}}_V$};
\end{tikzpicture}
\\[2em]
\begin{tikzpicture}[baseline=0, xscale=.4, yscale=.4]
\draw[thick,->] (0,2) to (0,0);
\node at (1.25,0) {$\phantom{^*}V^*$};
\node at (1.25,2) {$\phantom{^*}V^*$};
\draw[->] (1.25,.5) to (1.25,1.5);
\node[right] at (1.25,1) {$\mathrm{id}_{V^*}$};
\end{tikzpicture}
\qquad
\begin{tikzpicture}[baseline=0, xscale=.4, yscale=.4]
\draw[thick,->] (5,0) to [out=90, in=-90] (7,2);
\draw[line width=3pt, white] (7,0) to [out=90, in=-90] (5,2);
\draw[thick,->] (7,0) to [out=90, in=-90] (5,2);
\node[below] at (9.5,.5) {$\phantom{V}W\otimes V\phantom{W}$};
\node[above] at (9.5,1.5) {$\phantom{W}V\otimes W\phantom{V}$};
\draw[->] (9.5,.5) to (9.5,1.5);
\node[right] at (9.5,1) {$R_{VW}^{-1}$};
\end{tikzpicture}
\qquad
\begin{tikzpicture}[baseline=0, xscale=.4, yscale=.4]
\draw[thick,->] (6,2) to [out=-90, in=0] (5,.5);
\draw[thick,] (5,.5) to [out=180, in=-90] (4,2);
\node[above] at (8,1.5) {$\phantom{^*}V\otimes V^*$};
\node[below] at (8,.5) {$\mathbb{K}$};
\draw[->] (8,.5) to (8,1.5);
\node[right] at (8,1) {$\overleftarrow{\operatorname{coev}}_V$};
\end{tikzpicture}
\qquad
\begin{tikzpicture}[baseline=0, xscale=.4, yscale=.4]
\draw[thick,->] (4,2) to [out=-90, in=180] (5,.5);
\draw[thick,] (5,.5) to [out=0, in=-90] (6,2);
\node[above] at (8,1.5) {$V^*\otimes V\phantom{^*}$};
\node[below] at (8,.5) {$\mathbb{K}$};
\draw[->] (8,.5) to (8,1.5);
\node[right] at (8,1) {$\overrightarrow{\operatorname{coev}}_V$};
\end{tikzpicture}
\caption{A graphical definition of the Reshetikhin--Turaev functor on oriented
 elementary tangle diagrams.}\label{fig:FIG3}
\end{figure}

Set $V$ a vector space over a field $\mathbb{K}$ of characteristic zero and consider
$V^\ast = \mathrm{End}(V,\mathbb{K})$ its dual. Each elementary oriented tangle is associated
to a linear map as shown in Figure~\ref{fig:FIG3} under the Reshetikhin--Turaev functor
$F$. Now composing these elementary maps, one obtains a linear map $F_T$
for any oriented tangle diagram $T$.
For example, if $T$ is an oriented $(n,m)$-tangle diagram, then
$
F_{T} \in \mathrm{End}(V^{\epsilon_1} \otimes V^{\epsilon_2} \otimes \dots
\otimes V^{\epsilon_n}, V^{\delta_1} \otimes V^{\delta_2} \otimes \dots
\otimes V^{\delta_m}) ,
$
where $\epsilon_i$, $\delta_j$ are empty or $\ast$, depending on the
orientations of the boundary $\partial T$ of $T$.

\subsection{Invariants of tangles}
\label{sub:tangles}

To state the identities of an enhanced $R$-matrix, we need to fix some notation.
We do so following Ohtsuki~\cite{Ohtsuki} rather than Kashaev~\cite{GK:multi, Kas1}.
Fix a finite-dimensional vector space $V$ with basis $(e_1, \dots, e_n)$.
For $A \in \mathrm{End}(V \otimes V)$, we write $A = (A_{(i,j)}^{(k,l)})$ as a matrix in
basis $(e_i \otimes e_j)$. Therefore,
$A_{(i,j)}^{(k,l)} = (e_k \otimes e_l)^{\ast}(A(e_i \otimes e_j))$. Then we define
matrices $A^{\circlearrowleft}$ and $A^{\circlearrowright}$ by setting:
\begin{equation}
\label{Ab}
(A^\circlearrowleft)_{(k,i)}^{(l,j)} =A_{(i,j)}^{(k,l)}, \qquad
(A^\circlearrowright)_{(j,l)}^{(i,k)}=A_{(i,j)}^{(k,l)} .
\end{equation}

The next definition plays a crucial role in our paper {and is originally due to Turaev \cite{Turaev88,Turaev89}. Also see \cite[Section~XII.4]{Kassel}}.

\begin{Definition}
\label{def.Rh}
An \emph{enhanced $R$-matrix} on a finite-dimensional vector space $V$
is a pair of invertible endomorphisms $R \in \mathrm{End}(V\otimes V)$,
the \emph{$R$-matrix}, and $h \in \mathrm{End}(V)$, the \emph{enhancement}, satisfying
\begin{subequations}\label{eq1}
\begin{gather}
\label{eq1a}
R \circ (h \otimes h) = (h \otimes h) \circ R, \\
\label{eq1b}
\mathrm{tr}_{2} \bigl((\mathrm{id}_V \otimes h) \circ R^{\pm 1}\bigr) = \mathrm{id}_V, \\
\label{eq1c}
\bigl(R^{-1}\bigr)^{\circlearrowleft} \circ \bigl((\mathrm{id}_V \otimes h)
\circ R \circ \bigl(h^{-1} \otimes \mathrm{id}_V\bigr)\bigr)^{\circlearrowright} =
\mathrm{id}_V \otimes \mathrm{id}_{V^{\ast}}, \\
\label{eq1d}
(R \otimes \mathrm{id}_V) \circ (\mathrm{id}_V \otimes R)
\circ (R \otimes \mathrm{id}_V) = (\mathrm{id}_V \otimes R)
\circ (R \otimes \mathrm{id}_V) \circ (\mathrm{id}_V \otimes R).
\end{gather}
\end{subequations}
\end{Definition}

The partial trace $\mathrm{tr}_2$ of $f \otimes g \in \mathrm{End}(V \otimes V)$
is defined by
\begin{equation}
\label{tr2}
\mathrm{tr}_{2}(f \otimes g) = \mathrm{tr}(g) f \in \mathrm{End}(V)
\end{equation}
with the natural identification of $\mathrm{End}(V \otimes V)$ and
$\mathrm{End}(V) \otimes \mathrm{End}(V)$. {This extends to partial trace operations
$\mathrm{tr}_i\colon \mathrm{End}(V^{\otimes n})\to \mathrm{End}\bigl(V^{\otimes n-1}\bigr)$ where we trace in the $i$-th
component. Let $\mathrm{tr}_{i_1i_2\cdots i_j}$ denote the composition of partial traces
$\mathrm{tr}_{i_1}\circ\mathrm{tr}_{i_2}\circ\cdots\circ\mathrm{tr}_{i_j}$.}

For $1\leq i <j \leq n$ and $f\in\mathrm{End}(V\otimes V)$, define
$(f)_{ij}\in\mathrm{End}(V^{\otimes n})$ which acts by $f$ in the $i$-th and $j$-th
tensor factors and is the identity otherwise.

Given an $R$-matrix, equation \eqref{eq1d} defines a representation $\rho_R$ of
the braid group $B_n$ by mapping the elementary braid
generator in position $(k, k+1)$ to $(R)_{k,k+1}\in \mathrm{End}(V^{\otimes n})$.

\begin{Remark}
\label{rem.cups}
An enhanced $R$-matrix determines an operator valued invariant of isotopy classes
of tangles under the Reshetikhin--Turaev functor under the following definitions of
cup and cap maps:
\begin{alignat*}{3}
& \overrightarrow{\operatorname{ev}}_V(x \otimes f) := f(h(x)), \qquad &&
\overleftarrow{\operatorname{ev}}_V(f \otimes x) := f(x), &\\
& \overleftarrow{\operatorname{coev}}_V(1) := \sum_{i} e_i \otimes e_i^\ast, \qquad &&
\overrightarrow{\operatorname{coev}}_V(1) := \sum_{i} e_i^\ast \otimes h^{-1}(e_i)^\ast .&
\end{alignat*}
\end{Remark}

\begin{Remark}
\label{rem.Lcut}
If $\beta\in B_n$ is a braid with closure $L$ {and $V$ is simple}, then
\begin{gather*}
F_{L^{\mathrm{cut}}} =
\mathrm{tr}_{2,\dots, n}\bigl(\bigl(\mathrm{id}_V\otimes
h^{\otimes (n-1)}\bigr)\circ \rho_R(\beta) \bigr)\in \mathrm{End}(V),
\\
\langle F_{L^{\mathrm{cut}}}\rangle
=
\frac{1}{\dim(V)}
\mathrm{tr}\bigl(\bigl(\mathrm{id}_V\otimes
h^{\otimes (n-1)}\bigr)\circ \rho_R(\beta)\bigr) .
\end{gather*}
\end{Remark}

\subsection{Rotated tangles}

In this subsection we compare the tangles of the enhanced $R$-matrices
(see Figure~\ref{fig:tangles1}) with the rotated tangles used in~\cite{GK:multi}. The latter
are compositions of four types of segments
\begin{equation}
\label{4arcs}
\begin{tikzpicture}[yscale=.4,baseline]
\draw[thick,->] (0,0) to [out=90,in=-90] (0,1);
\end{tikzpicture}\ ,\
\begin{tikzpicture}[yscale=.4,baseline]
\draw[thick,<-] (0,0) to [out=90,in=-90] (0,1);
\end{tikzpicture}\ ,\ \quad
\begin{tikzpicture}[xscale=.7,yscale=1,baseline]
\draw[thick,<-] (0,0) to [out=90,in=90] (1,0);
\end{tikzpicture}\ ,\
\begin{tikzpicture}[xscale=.7,yscale=1,baseline=20]
\draw[thick,<-] (0,1) to [out=-90,in=-90] (1,1);
\end{tikzpicture}
\end{equation}
and eight types of crossings (four positive and four negative ones)
\begin{equation}
\label{8x}
\begin{tikzpicture}[scale=.7,baseline]
\draw[thick,<-] (0,1) to [out=-90,in=90] (1,0);
\draw[line width=3pt,white] (1,1) to [out=-90,in=90] (0,0);
\draw[thick,<-] (1,1) to [out=-90,in=90] (0,0);
\end{tikzpicture}\ ,\
\begin{tikzpicture}[scale=.7,baseline]
\draw[thick,<-] (1,1) to [out=-90,in=90] (0,0);
\draw[line width=3pt,white] (0,1) to [out=-90,in=90] (1,0);
\draw[thick,<-] (0,1) to [out=-90,in=90] (1,0);
\end{tikzpicture}\ ,\
\quad
\begin{tikzpicture}[scale=.7,baseline]
\draw[thick,->] (1,1) to [out=-90,in=90] (0,0);
\draw[line width=3pt,white] (0,1) to [out=-90,in=90] (1,0);
\draw[thick,<-] (0,1) to [out=-90,in=90] (1,0);
\end{tikzpicture}\ ,\
\begin{tikzpicture}[scale=.7,baseline]
\draw[thick,<-] (0,1) to [out=-90,in=90] (1,0);
\draw[line width=3pt,white] (1,1) to [out=-90,in=90] (0,0);
\draw[thick,->] (1,1) to [out=-90,in=90] (0,0);
\end{tikzpicture}\ ,\
\quad
\begin{tikzpicture}[scale=.7,baseline]
\draw[thick,->] (0,1) to [out=-90,in=90] (1,0);
\draw[line width=3pt,white] (1,1) to [out=-90,in=90] (0,0);
\draw[thick,->] (1,1) to [out=-90,in=90] (0,0);
\end{tikzpicture}\ ,\
\begin{tikzpicture}[scale=.7,baseline]
\draw[thick,->] (1,1) to [out=-90,in=90] (0,0);
\draw[line width=3pt,white] (0,1) to [out=-90,in=90] (1,0);
\draw[thick,->] (0,1) to [out=-90,in=90] (1,0);
\end{tikzpicture}\ ,\
\quad
\begin{tikzpicture}[scale=.7,baseline]
\draw[thick,<-] (1,1) to [out=-90,in=90] (0,0);
\draw[line width=3pt,white] (0,1) to [out=-90,in=90] (1,0);
\draw[thick,->] (0,1) to [out=-90,in=90] (1,0);
\end{tikzpicture}\ ,\
\begin{tikzpicture}[scale=.7,baseline]
\draw[thick,->] (0,1) to [out=-90,in=90] (1,0);
\draw[line width=3pt,white] (1,1) to [out=-90,in=90] (0,0);
\draw[thick,<-] (1,1) to [out=-90,in=90] (0,0);
\end{tikzpicture} .
\end{equation}

In addition, the remaining two types of segments are allowed, but when they
occur, they are replaced as follows:
\begin{equation}
\label{2racs}
\
\begin{tikzpicture}[baseline=5, xscale=0.4, yscale=0.45]
\draw[thick] (0,0) to [out=90,in=180] (1,1);
\draw[thick,->] (1,1) to [out=0,in=90] (2,0);
\end{tikzpicture}
\ \mapsto
\begin{tikzpicture}[baseline=5,xscale=.5,yscale=0.15]
\coordinate (a0) at (0,0);
\coordinate (a1) at (1,3);
\coordinate (a2) at (2,0);
\draw[thick] (a1) to [out=180,in=135] (a2);
\draw[line width=3, color=white] (a0) to [out=45,in=0] (a1);
\draw[thick,->] (a0) to [out=45,in=0] (a1);
\end{tikzpicture}
 , \qquad\qquad
\begin{tikzpicture}[baseline=-7, xscale=0.4, yscale=0.45]
\draw[thick] (0,0) to [out=-90,in=180] (1,-1);
\draw[thick,->] (1,-1) to [out=0,in=-90] (2,0);
\end{tikzpicture}\
\mapsto
\begin{tikzpicture}[baseline=-7,xscale=.5,yscale=0.15]
\coordinate (a0) at (0,0);
\coordinate (a1) at (1,-3);
\coordinate (a2) at (2,0);
\draw[thick] (a1) to [out=180,in=-135] (a2);
\draw[line width=3, color=white] (a0) to [out=-45,in=0] (a1);
\draw[thick,->] (a0) to [out=-45,in=0] (a1);
\end{tikzpicture} .
\end{equation}

The local weights are given by
\begin{equation*}
\begin{tikzpicture}[yscale=.7,xscale=.4,baseline=7]
\draw[thick,<-] (0,1) to [out=-90,in=90] (1,0);
\draw[line width=3pt,white] (1,1) to [out=-90,in=90] (0,0);
\draw[thick,<-] (1,1) to [out=-90,in=90] (0,0);
\node (sw) at (0,-.2){\tiny $a$};\node (se) at (1,-.2){\tiny $b$};
\node (nw) at (0,1.2){\tiny $c$};\node (ne) at (1,1.2){\tiny $d$};
\end{tikzpicture},
\begin{tikzpicture}[yscale=.7,xscale=.4,baseline=7]
\draw[thick,->] (1,1) to [out=-90,in=90] (0,0);
\draw[line width=3pt,white] (0,1) to [out=-90,in=90] (1,0);
\draw[thick,<-] (0,1) to [out=-90,in=90] (1,0);
\node (sw) at (0,-.2){\tiny $c$};\node (se) at (1,-.2){\tiny $a$};
\node (nw) at (0,1.2){\tiny $d$};\node (ne) at (1,1.2){\tiny $b$};
\end{tikzpicture},
\begin{tikzpicture}[yscale=.7,xscale=.4,baseline=7]
\draw[thick,->] (0,1) to [out=-90,in=90] (1,0);
\draw[line width=3pt,white] (1,1) to [out=-90,in=90] (0,0);
\draw[thick,->] (1,1) to [out=-90,in=90] (0,0);
\node (sw) at (0,-.2){\tiny $d$};\node (se) at (1,-.2){\tiny $c$};
\node (nw) at (0,1.2){\tiny $b$};\node (ne) at (1,1.2){\tiny $a$};
\end{tikzpicture}
\xmapsto{w_s}
\langle c^*\otimes d^*,R(a\otimes b)\rangle,
\qquad
\begin{tikzpicture}[yscale=.7,xscale=.4,baseline=7]
\draw[thick,<-] (1,1) to [out=-90,in=90] (0,0);
\draw[line width=3pt,white] (0,1) to [out=-90,in=90] (1,0);
\draw[thick,->] (0,1) to [out=-90,in=90] (1,0);
\node (sw) at (0,-.2){\tiny $b$};\node (se) at (1,-.2){\tiny $d$};
\node (nw) at (0,1.2){\tiny $a$};\node (ne) at (1,1.2){\tiny $c$};
\end{tikzpicture}\
\xmapsto{w_s}
\big\langle a\otimes c^*,{\bigl(\bigl(R^{-1}\bigr)^{\circlearrowleft}\bigr)^{-1}}(b\otimes d^*)\big\rangle
\end{equation*}
for positive crossings and likewise for negative crossings.
By \eqref{eq1c}, there is an equality
\[
\bigl(\bigl(R^{-1}\bigr)^{\circlearrowleft}\bigr)^{-1}=\bigl((\mathrm{id}_V \otimes h)
\circ R \circ \bigl(h^{-1} \otimes \mathrm{id}_V\bigr)\bigr)^{\circlearrowright}.
\]

The reason for using the rotated tangles of equations~\eqref{4arcs}, \eqref{8x}
and \eqref{2racs} rather than the ones in Figure~\ref{fig:tangles1} is that this
leads to an efficient computation of the corresponding link invariants using
the local methods and thin position highlighted by Bar-Natan and van der
Veen~\cite{BNV:PG,BNV:API} and implemented by Li and one of the authors
in~\cite{GL:patterns}.

Using the rigid $R$-matrix and the rotation numbers $\varphi=(1,-1,-1,1)$ in the
\texttt{Rot.m} program of~\cite[Section~2]{BNV:PG}, we can confirm by an explicit
calculation that the polynomial equations (11)--(13) of~\cite{BNV:PG} are satisfied
and hence one obtains invariants of tangles (up to isotopy) with no closed
components, hence endomorphism-valued invariants of (long) knots. By their very
definition, these invariants coincide with the ones of~\cite{GK:multi} since the weights
of the rotated crossings match in both cases.

Finally, to identify the long knot polynomial invariants of~\cite{GK:multi} with those
of Definition~\ref{def.Rh}, we use the fact that the~\cite{GK:multi}-weights of the
crossings of Figure~\ref{8x} are obtained by rotating the two tangles in the
right of Figure~\ref{fig:tangles1} accordingly and multiplying by the corresponding
caps/cups.

\subsection[The Alexander R-matrix]{The Alexander $\boldsymbol{R}$-matrix}
\label{sub:alex}

In this subsection, we review in detail the Alexander $R$-matrix, its enhancement,
its invariance under~\eqref{item:P1} and~\eqref{item:P2} and the identification
of the corresponding link invariant with the Alexander polynomial, following
the beautiful presentation of Ohtsuki~\cite{Ohtsuki}.

We begin by defining the enhanced $R$-matrix we use to define the Alexander polynomial.
Let~$Y$ denote the 2-dimensional $\mathbb{C}$-vector space with basis
$(v_0, v_1)$. Then $Y \otimes Y$ is equipped with the basis
$(v_0 \otimes v_0, v_0 \otimes v_1, v_1 \otimes v_0, v_1 \otimes v_1)$.
Define
\begin{equation}
\label{alexRh}
R_t = t^{-1/2}~\begin{pmatrix}
1 & 0 & 0 & 0 \\
0 & 0 & 1 & 0 \\
0 & t & 1-t & 0 \\
0 & 0 & 0 & -t
\end{pmatrix} \in \mathrm{End}(Y \otimes Y), \qquad h_t = t^{1/2}~\begin{pmatrix}
1 & 0 \\
0 & -1
\end{pmatrix} \in \mathrm{End}(Y) .
\end{equation}

The next lemma can be proven by a straightforward matrix computation.

\begin{Lemma}
\label{lem.alex0}
The pair $(R_t,h_t)$ satisfies equations~\eqref{eq1}.
\end{Lemma}

As discussed in the introduction, this leads to a definition of an operator-valued
invariant of isotopy class of tangles~$T$, which for reasons that will
become clear later, we will denote by~$F_{\Delta,T}$.

\begin{Lemma}
\label{lem.alex1}
The enhanced Alexander $R$-matrix satisfies properties \eqref{item:P1} and
\eqref{item:P2}.
\end{Lemma}

\begin{proof}
The proof uses a degree argument for the $R$-matrix (and its inverse)
and the tangle isotopies of Figure~\ref{fig:egalite} (for \eqref{item:P1})
and~\ref{fig:egalite2} (for \eqref{item:P2}).
We use the grading of the basis of $Y$ {and its dual}
$\mathrm{deg}(v_0) = 0$, $\mathrm{deg}(v_1) = 1$, $\mathrm{deg}(v_0^\ast) = 0$, $\mathrm{deg}(v_1^\ast) = -1$.
It follows by inspection that the $R$-matrix
and its inverse (as well as the cups/caps {with the convention that the ground field $\mathbb{C}$ is in degree $0$}) are degree-preserving.
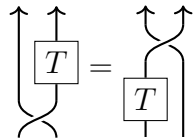
\begin{figure}[htpb!]
\centering\begin{tikzpicture}[yscale=1,xscale=1,baseline=-5]
\node[draw] (a) at (0,0) {$T$};
\draw[thick,->] (0,-.5)--(a)--(0,.75);
\draw[thick] (0,-1) to [out=90,in=-90] (-.5,-.5);
\draw[line width=3pt,white] (0,-.5) to [out=-90,in=90] (-.5,-1);
\draw[thick] (0,-.5) to [out=-90,in=90] (-.5,-1);
\draw[thick,->] (-.5,-.5) to (-.5,.75);
\end{tikzpicture}
=
\begin{tikzpicture}[yscale=1,xscale=1,baseline=9]
\node[draw] (a) at (0,0) {$T$};
\draw[thick] (0,-.5)--(a)--(0,.5);
\draw[thick] (.5,.5) to [out=90,in=-90] (0,1);
\draw[line width=3pt,white] (0,.5) to [out=90,in=-90] (.5,1);
\draw[thick] (0,.5) to [out=90,in=-90] (.5,1);
\draw[thick] (.5,-.5) to [out=90,in=-90] (.5,.5);
\draw[thick,->] (.5,1)--(.5,1.25);
\draw[thick,->] (0,1)--(0,1.25);
\end{tikzpicture}
\caption{A tangle isotopy.}
\label{fig:egalite}
\end{figure}

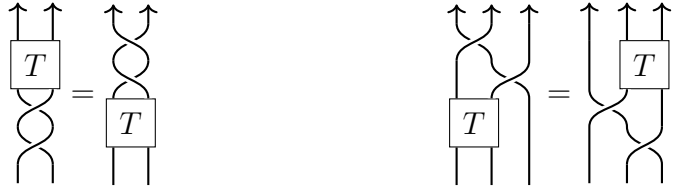
\begin{figure}[htpb!]
\centering
\begin{tikzpicture}[yscale=1,xscale=1,baseline=0]
\node[rectangle,inner sep=5pt,draw] (c) at (1/2,1/2) {$T$};
\draw[thick,<-] ($(c.north east)+(-.1,1/2)$)--($(c.north east)+(-.1,0)$);
\draw[thick,<-] ($(c.north west)+(.1,1/2)$)--($(c.north west)+(.1,0)$);
\draw[thick] ($(c.south west)+(.1,0)$) to [out=-90,in=90]
($(c.south east)+(-.1,-1/2)$);
\draw[line width=3pt,white] ($(c.south east)+(-.1,0)$) to [out=-90,in=90]
($(c.south west)+(.1,-1/2)$);
\draw[thick] ($(c.south east)+(-.1,0)$) to [out=-90,in=90]
($(c.south west)+(.1,-1/2)$);
\draw[thick] ($(c.south west)+(.1,-1/2)$) to [out=-90,in=90]
($(c.south east)+(-.1,-1)$);
\draw[line width=3pt,white] ($(c.south east)+(-.1,-1/2)$) to [out=-90,in=90]
($(c.south west)+(.1,-1)$);
\draw[thick] ($(c.south east)+(-.1,-1/2)$) to [out=-90,in=90]
($(c.south west)+(.1,-1)$);
\draw[thick] ($(c.south east)+(-.1,-1)$)--($(c.south east)+(-.1,-1)+(0,-.25)$);
\draw[thick] ($(c.south west)+(.1,-1)$)--($(c.south west)+(.1,-1)+(0,-.25)$);
\end{tikzpicture}
=
\begin{tikzpicture}[yscale=1,xscale=1,baseline=22]
\node[rectangle,inner sep=5pt,draw] (c) at (1/2,1/2) {$T$};
\draw[thick,] ($(c.south east)+(-.1,-1/2)$)--($(c.south east)+(-.1,0)$);
\draw[thick] ($(c.south west)+(.1,-1/2)$)--($(c.south west)+(.1,0)$);
\draw[thick] ($(c.north west)+(.1,0)$) to [out=90,in=-90]
($(c.north east)+(-.1,1/2)$);
\draw[line width=3pt,white] ($(c.north east)+(-.1,0)$)
to [out=90,in=-90] ($(c.north west)+(.1,1/2)$);
\draw[thick] ($(c.north east)+(-.1,0)$) to [out=90,in=-90]
($(c.north west)+(.1,1/2)$);
\draw[thick] ($(c.north west)+(.1,1/2)$) to [out=90,in=-90]
($(c.north east)+(-.1,1)$);
\draw[line width=3pt,white] ($(c.north east)+(-.1,1/2)$)
to [out=90,in=-90] ($(c.north west)+(.1,1)$);
\draw[thick] ($(c.north east)+(-.1,1/2)$) to [out=90,in=-90]
($(c.north west)+(.1,1)$);
\draw[thick,->] ($(c.north east)+(-.1,1)$) -- ($(c.north east)+(-.1,1)+(0,.25)$);
\draw[thick,->] ($(c.north west)+(.1,1)$) -- ($(c.north west)+(.1,1)+(0,.25)$);
\end{tikzpicture}
\quad\quad
\begin{tikzpicture}[yscale=1,xscale=1,baseline=22]
\node[rectangle,inner sep=5pt,draw] (c) at (1/2,1/2) {$T$};
\draw[thick] ($(c.south east)+(-.1,-1/2)$)--($(c.south east)+(-.1,0)$);
\draw[thick] ($(c.south west)+(.1,-1/2)$)--($(c.south west)+(.1,0)$);
\draw[thick] ($(c.south east)+(.4,-1/2)$) to [out=90,in=-90]
($(c.north east)+(.4,0)$) to [out=90,in=-90] ($(c.north east)+(-.1,1/2)$)
to [out=90,in=-90] ($(c.north west)+(.1,1)$);
\draw[thick,->] ($(c.north west)+(.1,1)$) to ($(c.north west)+(.1,1.25)$);
\draw[line width=3pt,white] ($(c.north east)+(-.1,0)$) to [out=90,in=-90]
($(c.north east)+(.4,1/2)$);
\draw[thick] ($(c.north east)+(-.1,0)$) to [out=90,in=-90]
($(c.north east)+(.4,1/2)$) to [out=90,in=-90] ($(c.north east)+(.4,1)$);
\draw[thick,->] ($(c.north east)+(.4,1)$) to ($(c.north east)+(.4,1.25)$);
\draw[thick] ($(c.north west)+(.1,0)$) to [out=90,in=-90]
($(c.north west)+(.1,1/2)$);
\draw[line width=3pt,white] ($(c.north west)+(.1,1/2)$) to [out=90,in=-90]
($(c.north east)+(-.1,1)$);
\draw[thick] ($(c.north west)+(.1,1/2)$) to [out=90,in=-90]
($(c.north east)+(-.1,1)$);
\draw[thick,->] ($(c.north east)+(-.1,1)$) to ($(c.north east)+(-.1,1.25)$);
\end{tikzpicture}
=
\begin{tikzpicture}
[yscale=1,xscale=1,baseline=0]
\node[rectangle,inner sep=5pt,draw] (c) at (1/2,1/2) {$T$};
\draw[thick,<-] ($(c.north east)+(-.1,1/2)$)--($(c.north east)+(-.1,0)$);
\draw[thick,<-] ($(c.north west)+(.1,1/2)$)--($(c.north west)+(.1,0)$);
\draw[thick] ($(c.south east)+(-.1,-1.25)$) to [out=90,in=-90]
($(c.south east)+(-.1,-1.)$) to [out=90,in=-90] ($(c.south west)+(.1,-1/2)$)
to [out=90,in=-90] ($(c.south west)+(-.4,0)$) to [out=90,in=-90,]
($(c.north west)+(-.4,0)$);
\draw[thick,->] ($(c.north west)+(-.4,0)$)--($(c.north west)+(-.4,1/2)$);
\draw[thick] ($(c.south east)+(-.1,0)$) to [out=-90,in=90]
($(c.south east)+(-.1,-1/2)$);
\draw[line width=3pt,white] ($(c.south east)+(-.1,-1/2)$) to [out=-90,in=90]
($(c.south west)+(.1,-1)$);
\draw[thick] ($(c.south east)+(-.1,-1/2)$) to [out=-90,in=90]
($(c.south west)+(.1,-1)$);
\draw[line width=3pt,white] ($(c.south west)+(.1,0)$) to [out=-90,in=90]
($(c.south west)+(-.4,-1/2)$);
\draw[thick] ($(c.south west)+(.1,0)$) to [out=-90,in=90]
($(c.south west)+(-.4,-1/2)$) to ($(c.south west)+(-.4,-1.25)$) ;
\draw[thick] ($(c.south west)+(.1,-1)$) to ($(c.south west)+(.1,-1.25)$);
\end{tikzpicture}
\caption{More tangle isotopies.}
\label{fig:egalite2}
\end{figure}

Let us first prove that $R_{t}$ satisfies \eqref{item:P1}. Fix $T$ a
$(1,1)$-tangle. Using the degree-preservation of $F_{\Delta,T} \in \mathrm{End}(Y)$,
it follows that
\begin{equation}
\label{T1Alex}
F_{\Delta,T} = \begin{pmatrix}
\alpha & 0 \\
0 & \beta
\end{pmatrix} .
\end{equation}

The invariance of $F_\Delta$ under the isotopy of
Figure~\ref{fig:egalite} implies that
\begin{equation*}
(\mathrm{id}_Y \otimes F_{\Delta,T}) \circ R_{t} =
R_{t} \circ (F_{\Delta,T} \otimes \mathrm{id}_Y) .
\end{equation*}
Inserting~\eqref{T1Alex} in the above identity shows that $\alpha = \beta$. Therefore,
$ F_{\Delta,T} = \alpha\,\mathrm{id}_Y$. So $R_{t}$ satisfies~\eqref{item:P1}.

Now we prove that $R_{t}$ satisfies \eqref{item:P2}. Fix a $(2,2)$-tangle $T$
and consider its left and right closures $T_1$ and $T_2$, as in
Figure~\ref{fig:T12}. Using the degree-preservation of the $R$-matrix, it follows
that the operator invariant representing $T$ has the following form when
written in the standard basis for $Y \otimes Y$:
\[
F_{\Delta,T}=~\begin{pmatrix}
a & 0 & 0 & 0 \\
0 & b & c & 0 \\
0 & d & e & 0 \\
0 & 0 & 0 & f
\end{pmatrix},
\]
where $a$, $b$, $c$, $d$, $e$, $f$ are $6$ parameters. The tangle isotopy of
the left- and right-hand sides of Figure~\ref{fig:egalite2} implies the equations
\begin{gather*}
 R_t^2 \circ F_{\Delta,T} - F_{\Delta,T} \circ R_t^2 = 0 ,\\
 (R_t \otimes \text{id}_Y) \circ (\text{id}_Y \otimes R_t) \circ
( F_{\Delta,T} \otimes \text{id}_Y) - (\text{id}_Y \otimes F_{\Delta,T})
\circ (R_t \otimes \text{id}_Y) \circ (\text{id}_Y \otimes R_t) = 0 .
\end{gather*}
This is a linear system of equations in $6$ variables $a,b,\dots,f$ with
coefficients in $\mathbb{Q}(t)$. It has rank $2$ and solving, we find that
$a = b + c$, $d = c t$, $e = b + c (1-t)$, $f = b - c t$.
Using this and the definition of $\mathrm{tr}_2$ from equation~\eqref{tr2}, we can now
compute $F_{\Delta,T_1}$ and $F_{\Delta,T_2}$, and we find
\begin{gather*}
\begin{split}
& F_{\Delta,T_1} =
 \mathrm{tr}_{2} ((\mathrm{id}_Y \otimes h_t) \circ F_{\Delta,T})
= c t^{1/2}\,\mathrm{id}_Y, \\
& F_{\Delta,T_2} = \mathrm{tr}_{1} \bigl(F_{\Delta,T}
\circ \bigl(h_t^{-1} \otimes \mathrm{id}_Y\bigr)\bigr) = c t^{1/2}\,\mathrm{id}_Y.
\end{split}
\end{gather*}
It follows that $F_{\Delta,T_1}=F_{\Delta,T_2}$ and this
concludes the proof of \eqref{item:P2}.
\end{proof}

\begin{Lemma}\label{lem.alex2}
The corresponding link invariant from Lemma~{\rm \ref{lem.alex1}} is the
Alexander polynomial.
\end{Lemma}

\begin{proof}
Recall that the Alexander polynomial of links satisfies the
skein relation and the initial condition
\begin{equation*}
\Delta_{L_+}(t) - \Delta_{L_-}(t) = \bigl(t^{-1/2} - t^{1/2}\bigr) \Delta_{L_0}(t),
\qquad \Delta_\bigcirc(t) = 1,
\end{equation*}
where $(L_+,L_-,L_0)$ is a triple of links as in Figure~\ref{fig:skein}.

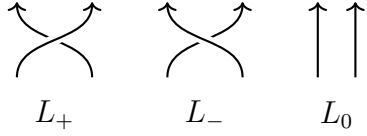
\begin{figure}[htpb!]\centering
\begin{tikzpicture}
[yscale=1,xscale=1,baseline=0]
\draw[thick,->] (1,0) to [out=90, in=-90] (0,1);
\draw[line width=4pt,white] (0,0) to [out=90, in=-90] (1,1);
\draw[thick,->] (0,0) to [out=90, in=-90] (1,1);
\node[] (a) at (1/2,-1/2) {$L_+$};
\draw[thick,->] (2,0) to [out=90, in=-90] (3,1);
\draw[line width=4pt,white] (3,0) to [out=90, in=-90] (2,1);
\draw[thick,->] (3,0) to [out=90, in=-90] (2,1);
\node[] (b) at (5/2,-1/2) {$L_-$};
\draw[thick,->] (4,0) to [out=90, in=-90] (4,1);
\draw[thick,->] (4.5,0) to [out=90, in=-90] (4.5,1);
\node[] (c) at (4.25,-1/2) {$L_0$};
\end{tikzpicture}
\caption{A triple of links.}
\label{fig:skein}
\end{figure}

It is well-known that the above skein theory uniquely characterizes the
Alexander polynomial.

The identity
$R_t - R_t^{-1} = \bigl(t^{-1/2} - t^{1/2}\bigr) \mathrm{id}_{Y \otimes Y}$
implies that the link invariant satisfies the same skein relation as
the Alexander polynomial, with the same initial condition. The result follows.
\end{proof}

\subsection[Conjugation of R-matrices]{Conjugation of $\boldsymbol{R}$-matrices}
\label{sub.conj}

In this section we make remarks concerning different $R$-matrices
leading to the same link invariants.
We say that two $R$-matrices $R \in \mathrm{End}\bigl(V^{\otimes 2}\bigr)$ and
$R' \in \mathrm{End}\bigl(W^{\otimes 2}\bigr)$ are \emph{conjugate} if there exists an isomorphism
$\varphi\colon V \to W$ such that
$
(\varphi \otimes \varphi) R \bigl(\varphi^{-1} \otimes \varphi^{-1}\bigr) = R'$.
Note that conjugate $R$-matrices on $V$ are conjugate automorphisms of $V \otimes V$,
but the converse is not true.
This definition extends naturally to enhanced $R$-matrices $(R,h)$ and $(R',h')$
by requiring $\phi h \phi^{-1} = h'$. If one $R$-matrix is enhanced and
conjugate to another $R$-matrix, then that one is canonically enhanced too.

It is easy to see that the operator-valued RT invariants of two conjugate
enhanced $R$-matrices are conjugate, and hence if one satisfies \eqref{item:P1}
and \eqref{item:P2},
so does the other, and the two associated link invariants are equal.
But equality of link invariants is also implied by a weaker version of conjugacy
for $R$-matrices, a fact that should be better-known.

We say that two enhanced $R$-matrices $(R,h)$ and $(R',h')$, where
$R \in \mathrm{End}\bigl(V^{\otimes 2}\bigr)$ and $R' \in \mathrm{End}\bigl(W^{\otimes 2}\bigr)$, are
\emph{weakly-conjugate} if there exist automorphisms
$\varphi_n\colon V^{\otimes n}\to W^{\otimes n}$ for all $n\geq 1$ such that
\begin{enumerate}\itemsep=0pt \setlength{\leftskip}{0.26cm}
\item[(\namedlabel{item:BC1}{$BC_1$})]
 The maps determine an equivalence of braid group representations, i.e.,
$\varphi_n \rho_R(\beta) \varphi_n^{-1}= \rho_{R'}(\beta)$.
\item[(\namedlabel{item:BC2}{$BC_2$})]
 There exist automorphisms $\sigma \in \mathrm{Hom}(V,W)$, $\nu_{n},\gamma_n\in \mathrm{End}(V^{\otimes n})$
 which satisfy the following properties for $n\geq 2$:
\begin{gather}
 \varphi_n=(\varphi_{n-1}\otimes \mathrm{id}_W)\circ (\nu_{n-1} \otimes \sigma)\circ\gamma_n,
 \label{eq.phin}
 \\
 \sigma^{-1} \circ h' \circ\sigma = h, \label{eq.sn}
 \\
 \mathrm{tr}_n\bigl((\mathrm{id}_{V^{\otimes (n-1)}}\otimes h)\circ \gamma_n\circ F_T\circ \gamma_n^{-1}\bigr)
 =\mathrm{tr}_n((\mathrm{id}_{V^{\otimes (n-1)}}\otimes h)\circ F_T),\label{eq.gn}
 \\
 \nu_{n-1} \circ \mathrm{tr}_n((\mathrm{id}_{V^{\otimes (n-1)}}\otimes h)\circ F_T)\circ \nu_{n-1}^{-1}
 =\mathrm{tr}_n((\mathrm{id}_{V^{\otimes (n-1)}}\otimes h)\circ F_T)\label{eq.nun}
\end{gather}
for any $(n,n)$-tangle $T$ with upward oriented boundary. When $n=2$, we also require
\begin{gather}
 \mathrm{tr}_1\bigl(\bigl(h^{-1}\otimes\mathrm{id}_V\bigr)\circ \gamma_2\circ F_T\circ \gamma_2^{-1}\bigr)
 =
 \mathrm{tr}_1\bigl(\bigl(h^{-1}\otimes\mathrm{id}_V\bigr)\circ F_T\bigr), \label{eq.g2}
 \\
 \nu_1^{-1} \circ h \circ \nu_1=h .\label{eq.n1}
\end{gather}
\end{enumerate}

\begin{Lemma}\label{lem.RTequal}
Suppose $(R,h)$ and $(R',h')$ are weakly-conjugate enhanced $R$-matrices. If the RT
invariant of links associated to $(R,h)$ determines a modified RT invariant, then
$(R',h')$ also determines a modified RT invariant and these two invariants are equal.
\end{Lemma}

\begin{proof}
It is sufficient to show that the RT invariant associated to $(R',h')$ satisfies
properties~\eqref{item:P1} and~\eqref{item:P2}. Let $L$ be a link with braid
representative $\beta\in B_n$ and its modified invariant $\text{mRT}_{L,R}$ associated to
$(R,h)$. The key observation is that the partial trace on a weak conjugacy implies
a weak conjugacy between partial traces. Then
\begin{equation}
\begin{tikzpicture}[yscale=1,xscale=1,baseline=12] 
\node[rectangle,draw] (c) at (1/2,1/2) {$\rho_{R'}(\beta)$};
\node at ($(c.south west)+(.55,-.3)$) {$\cdots$};
\node at ($(c.north west)+(.55,.2)$) {$\cdots$};
\draw[thick] ($(c.south east)+(-.5,-1/2)$)--($(c.south east)+(-.5,0)$);
\draw[thick,->] ($(c.north east)+(-.5,0)$)--($(c.north east)+(-.5,1/2)$);
\draw[thick,->] ($(c.north west)+(.1,0)$)--($(c.north west)+(.1,1/2)$);
\draw[thick] ($(c.south west)+(.1,0)$)--($(c.south west)+(.1,-1/2)$);
\draw[thick] ($(c.north east)+(-.1,0)$) to [out=90, in= 180]
($(c.north east)+(.4,1/2)$) to [out=0, in= 90] ($(c.north east)+(.9,0)$);
\draw[thick,->] ($(c.north east)+(.9,0)$) to ($(c.east)+(.9,0)$);
\draw[thick] ($(c.east)+(.9,0)$) to ($(c.south east)+(.9,0)$) to [out=-90, in= 0]
($(c.south east)+(.4,-1/2)$) to [out=180, in=-90] ($(c.south east)+(-.1,0)$);
\end{tikzpicture}
=
\begin{tikzpicture} 
[scale=1,baseline=12]
\node[rectangle,draw] (c) at (1/2,1/2) {$\phantom{\rho_{R'}(\beta)}$};
\node at (c) {$\rho_{R}(\beta)$};
\node[rectangle,draw] (i) at ($(c.south)+(0,-1/2)$) {$\phantom{\rho_{R'}(\beta)}$};
\node[rectangle,draw] (p) at ($(c.north)+(0,1/2)$) {$\phantom{\rho_{R'}(\beta)}$};
\node at (i) {$\varphi_n^{-1}$};
\node at (p) {$\varphi_n$};
\node at ($(i.south west)+(.55,-.3)$) {$\cdots$};
\node at ($(p.north west)+(.55,.2)$) {$\cdots$};
\node[] at ($(i.north west)+(.55,.05)$) {$\cdots$};
\node[] at ($(c.north west)+(.55,.05)$) {$\cdots$};
\draw[thick] ($(i.north east)+(-.1,0)$)--($(c.south east)+(-.1,0)$);
\draw[thick] ($(c.north east)+(-.1,0)$)--($(p.south east)+(-.1,0)$);
\draw[thick] ($(i.south east)+(-.5,-1/2)$)--($(i.south east)+(-.5,0)$);
\draw[thick] ($(i.south west)+(.1,-1/2)$)--($(i.south west)+(.1,0)$);
\draw[thick,->] ($(p.north east)+(-.5,0)$)--($(p.north east)+(-.5,1/2)$);
\draw[thick,->] ($(p.north west)+(.1,0)$)--($(p.north west)+(.1,1/2)$);
\draw[thick] ($(p.north east)+(-.1,0)$) to [out=90, in= 180]
($(p.north east)+(.4,1/2)$) to [out=0, in= 90] ($(p.north east)+(.9,0)$);
\draw[thick,->] ($(p.north east)+(.9,0)$) to ($(c.east)+(.9,0)$);
\draw[thick] ($(c.east)+(.9,0)$) to ($(i.south east)+(.9,0)$) to [out=-90, in= 0]
($(i.south east)+(.4,-1/2)$)
to [out=180, in= -90] ($(i.south east)+(-.1,0)$);
\draw[thick] ($(c.south east)+(-.5,0)$) to ($(i.north east)+(-.5,0)$);
\draw[thick] ($(c.north west)+(.1,0)$) to ($(p.south west)+(.1,0)$);
\draw[thick] ($(c.north east)+(-.5,0)$) to ($(p.south east)+(-.5,0)$);
\draw[thick] ($(i.north east)+(-.5,0)$) to ($(c.south east)+(-.5,0)$);
\draw[thick] ($(c.north east)+(-.5,0)$) to ($(p.south east)+(-.5,0)$);
\draw[thick] ($(c.south west)+(.1,0)$) to ($(i.north west)+(.1,0)$);
\end{tikzpicture}
=
\begin{tikzpicture} 
[yscale=1,xscale=1,baseline=12]
\node[rectangle,draw] (c) at (1/2,1/2) {\phantom{$\varphi_{n-1}^{-1}\rho_{R'}(\beta)$}};
\node[rectangle,draw] (g) at ($(c.north)+(0,1/2)$)
{\phantom{$\varphi_{n-1}^{-1}\rho_{R'}(\beta)$}};
\node[rectangle,draw] (gi) at ($(c.south)+(0,-1/2)$)
{\phantom{$\varphi_{n-1}^{-1}\rho_{R'}(\beta)$}};
\node[rectangle,draw,anchor = west] (n) at ($(g.north west)+(0,1/2)$)
{\phantom{$\varphi_{n-1}^{-1}$}};
\node[rectangle,draw,anchor = west] (ni) at ($(gi.south west)+(0,-1/2)$)
{\phantom{$\varphi_{n-1}^{-1}$}};
\node[rectangle,draw] (i) at ($(ni.south)+(0,-1/2)$) {\phantom{$\varphi_{n-1}^{-1}$}};
\node[rectangle,draw] (p) at ($(n.north)+(0,1/2)$) {\phantom{$\varphi_{n-1}^{-1}$}};
\node[rectangle,draw,anchor = east] (si) at ($(gi.south east)+(0,-1/2)$)
{\phantom{$\sigma^{-1}$}};
\node[rectangle,draw,anchor = east] (s) at ($(g.north east)+(0,1/2)$)
{\phantom{$\sigma^{-1}$}};

\node[] at (si) {$\sigma^{-1}$};
\node[] at (s) {$\sigma$};
\node[] at (i) {$\varphi_{n-1}^{-1}$};
\node[] at (p) {$\varphi_{n-1}$};
\node at ($(i.south west)+(.6,-.3)$) {$\cdots$};
\node at ($(i.north west)+(.6,.05)$) {$\cdots$};
\node at ($(ni.north west)+(.6,.05)$) {$\cdots$};
\node at ($(gi.north west)+(.6,.05)$) {$\cdots$};
\node at ($(p.north west)+(.6,.2)$) {$\cdots$};
\node at ($(p.south west)+(.6,-.07)$) {$\cdots$};
\node at ($(n.south west)+(.6,-.07)$) {$\cdots$};
\node at ($(g.south west)+(.6,-.07)$) {$\cdots$};
\draw[thick] ($(i.south west)+(.1,-1/2)$)--($(i.south west)+(.1,0)$);
\draw[thick] ($(i.south east)+(-.1,-1/2)$)--($(i.south east)+(-.1,0)$);
\draw[thick,->] ($(p.north east)+(-.1,0)$)--($(p.north east)+(-.1,1/2)$);
\draw[thick] ($(i.north east)+(-.1,0)$)--($(p.south east)+(-.1,0)$);
\draw[thick] ($(i.north west)+(.1,0)$)--($(p.south west)+(.1,0)$);
\draw[thick,->] ($(p.north west)+(.1,0)$)--($(p.north west)+(.1,1/2)$);
\draw[thick] ($(si.north)$) to ($(s.south)$);
\coordinate (R) at ($(s.north)+(1,0)$);
\draw[thick,->] ($(s.north)$) to [out=90, in=180]($(s.north)+(.5,1/2)$)
to [out=0, in=90] (R) to (R |- c);
\draw[thick] ($(si.south)$) to [out=-90, in=180]($(si.south)+(.5,-1/2)$)
to [out=0, in=-90] ($(si.south)+(1,0)$) to (R |- c);
\node[rectangle,draw,fill=white] at (c) {\phantom{$\varphi_{n-1}^{-1}\rho_{R'}(\beta)$}};
\node[rectangle,draw,fill=white] at (g) {\phantom{$\varphi_{n-1}^{-1}\rho_{R'}(\beta)$}};
\node[rectangle,draw,fill=white] at (gi) {\phantom{$\varphi_{n-1}^{-1}\rho_{R'}(\beta)$}};
\node[rectangle,draw,fill=white] at (n) {\phantom{$\varphi_{n-1}^{-1}$}};
\node[rectangle,draw,fill=white] at (ni) {\phantom{$\varphi_{n-1}^{-1}$}};
\node[] at (c) {$\rho_{R}(\beta)$};
\node[] at (g) {$\gamma_n$};
\node[] at (gi) {$\gamma_n^{-1}$};
\node[] at (n) {$\nu_{n-1}$};
\node[] at (ni) {$\nu_{n-1}^{-1}$};
\end{tikzpicture}
=
\begin{tikzpicture}
[yscale=1,xscale=1,baseline=12]
\node[rectangle,draw] (c) at (1/2,1/2) {\phantom{$\varphi_{n-1}^{-1}$\hspace{2ex}}};
\node[rectangle,draw,anchor = west] (i) at ($(c.south west)+(0,-1/2)$)
{\phantom{$\varphi_{n-1}^{-1}$}};
\node[] at (i) {$\varphi_{n-1}^{-1}$};
\node[rectangle,draw,anchor = west] (p) at ($(c.north west)+(0,1/2)$)
{\phantom{$\varphi_{n-1}^{-1}$}};
\node[] at (p) {$\varphi_{n-1}$};
\draw[thick,->] ($(i.south west)+(.1,-1/2)$)--($(p.north west)+(.1,1/2)$);
\draw[thick,->] ($(i.south east)+(-.1,-1/2)$)--($(p.north east)+(-.1,1/2)$);
\node at ($(i.south west)+(.6,-.3)$) {$\cdots$};
\node at ($(p.north west)+(.6,.2)$) {$\cdots$};
\node at ($(p.south west)+(.6,-.07)$) {$\cdots$};
\node at ($(i.north west)+(.6,.05)$) {$\cdots$};
\draw[thick,->] ($(c.north east)+(-.1,0)$) to [out=90, in= 180]
($(c.north east)+(.4,1/2)$) to [out=0, in= 90] ($(c.north east)+(.9,0)$)
to ($(c.east)+(.9,0)$);
\draw[thick] ($(c.east)+(.9,0)$) to ($(c.south east)+(.9,0)$)
to [out=-90, in= 0] ($(c.south east)+(.4,-1/2)$)
to [out=180, in= -90] ($(c.south east)+(-.1,0)$);
\node[rectangle,draw,fill=white] at (c) {\phantom{$\varphi_{n-1}^{-1}$\hspace{2ex}}};
\node[] at (c) {$\rho_{R}(\beta)$};
\node[rectangle,draw,fill=white] at (p) {\phantom{$\varphi_{n-1}^{-1}$}};
\node at (p) {$\varphi_{n-1}$};
\node[rectangle,draw,fill=white] at (i) {\phantom{$\varphi_{n-1}^{-1}$}};
\node at (i) {$\varphi_{n-1}^{-1}$};
\end{tikzpicture}
\end{equation}
The first equality is due to property \eqref{item:BC1}, the second equality
is due to equation~\eqref{eq.phin} of~\eqref{item:BC2}, and the last equality
follows from equations~\eqref{eq.sn}, \eqref{eq.gn}, and \eqref{eq.nun} of~\eqref{item:BC2}. Thus we have a graphical proof of property~\eqref{item:P1} by
induction
\begin{equation}
\label{eq.P1R'}
\begin{tikzpicture}[yscale=1,xscale=1,baseline=12]
\node[rectangle,draw] (c) at (1/2,1/2) {$\rho_{R'}(\beta)$};
\draw[thick] ($(c.south west)+(.1,-1.4)$) to ($(c.south west)+(.1,0)$);
\draw[thick,->] ($(c.north west)+(.1,0)$) to ($(c.north west)+(.1,1.4)$);
\draw[thick,->] ($(c.north east)+(-.1,0)$) to [out=90, in= 180]
($(c.north east)+(.4,1/2)$) to [out=0, in= 90] ($(c.north east)+(.9,0)$)
to ($(c.east)+(.9,0)$);
\draw[thick] ($(c.east)+(.9,0)$) to ($(c.south east)+(.9,0)$)
to [out=-90, in= 0] ($(c.south east)+(.4,-1/2)$)
to [out=180, in= -90] ($(c.south east)+(-.1,0)$);
\draw[thick,->] ($(c.north east)+(-1,0)$) to [out=90, in= 180]
($(c.north east)+(.4,1.4)$) to [out=0, in= 90] ($(c.north east)+(1.8,0)$)
to ($(c.east)+(1.8,0)$);
\draw[thick] ($(c.east)+(1.8,0)$) to ($(c.south east)+(1.8,0)$)
to [out=-90, in= 0] ($(c.south east)+(.4,-1.4)$)
to [out=180, in= -90] ($(c.south east)+(-1,0)$);
\node at ($(c.east)+(1.35,0)$) {$\cdots$};
\end{tikzpicture}
=
\begin{tikzpicture}[yscale=1,xscale=1,baseline=12]
\node[rectangle,draw] (c) at (1/2,1/2) {\phantom{$\rho_{R'}(\beta)$}};
\node at (c) {$\rho_{R}(\beta)$};
\node[rectangle,draw] (p) at ($(c.north west)+(.1,1.5)$) {\phantom{$\varphi_1^{-1}$}};
\node[rectangle,draw] (i) at ($(c.south west)+(.1,-1.5)$) {\phantom{$\varphi_1^{-1}$}};
\node at (p) {$\varphi_1$};
\node at (i) {$\varphi_1^{-1}$};
\draw[thick] ($(i.south)+(0,-1/2)$) to (i.south);\draw[thick] (i.north)
to ($(c.south west)+(.1,0)$);
\draw[thick] ($(c.north west)+(.1,0)$) to (p.south);
\draw[thick,->] (p.north) to ($(p.north)+(0,1/2)$);
\draw[thick,->] ($(c.north east)+(-.1,0)$) to [out=90, in= 180]
($(c.north east)+(.4,1/2)$) to [out=0, in= 90] ($(c.north east)+(.9,0)$)
to ($(c.east)+(.9,0)$);
\draw[thick] ($(c.east)+(.9,0)$) to ($(c.south east)+(.9,0)$)
to [out=-90, in= 0] ($(c.south east)+(.4,-1/2)$)
to [out=180, in= -90] ($(c.south east)+(-.1,0)$);
\draw[thick,->] ($(c.north east)+(-1,0)$) to [out=90, in= 180]
($(c.north east)+(.4,1.4)$) to [out=0, in= 90] ($(c.north east)+(1.8,0)$)
to ($(c.east)+(1.8,0)$);
\draw[thick] ($(c.east)+(1.8,0)$) to ($(c.south east)+(1.8,0)$)
to [out=-90, in= 0] ($(c.south east)+(.4,-1.4)$)
to [out=180, in= -90] ($(c.south east)+(-1,0)$);
\node at ($(c.east)+(1.35,0)$) {$\cdots$};
\end{tikzpicture}
=
\text{mRT}_{L,R}
\begin{tikzpicture}[yscale=1,xscale=1,baseline=12]
 \coordinate[] (c) at (1/2,1/2) ;
\node[rectangle,draw] (p) at ($(c)+(0,1/2)$) {\phantom{$\varphi_1^{-1}$}};
\node[rectangle,draw] (i) at ($(c)+(0,-1/2)$) {\phantom{$\varphi_1^{-1}$}};
\node at (p) {$\varphi_1$};
\node at (i) {$\varphi_1^{-1}$};
\draw[thick] ($(i.south)+(0,-1/2)$) to (i.south);\draw[thick] (i.north) to (p.south);
\draw[thick,->] (p.north) to ($(p.north)+(0,1/2)$);
\end{tikzpicture}
=
\text{mRT}_{L,R}
\begin{tikzpicture}[yscale=1,xscale=1,baseline=12]
\coordinate[] (c) at (1/2,1/2) ;
\draw[thick,] ($(c)+(0,-1)$) to (c);
\draw[thick,->] (c) to ($(c)+(0,1)$);
\end{tikzpicture}
\end{equation}
using property~\eqref{item:P1} of the RT invariant associated to $(R,h)$.
We now show property~\eqref{item:P2} for $(R',h')$ in the equalities below using
it for $(R,h)$ and equations \eqref{eq.g2}--\eqref{eq.P1R'}. Let $T$ be the closure of $\beta$ over the right $n-2$ stands with $F_{R',T}$ and $F_{R,T}$ its associated RT morphisms. Then
\begin{equation*}
\begin{aligned}
\begin{tikzpicture} 
[yscale=1,xscale=1,baseline=12]
\node[rectangle,draw] (c) at (1/2,1/2) {$F_{R',T}$};
\draw[thick] ($(c.south east)+(-.4,-.5)$) to ($(c.south east)+(-.4,0)$);
\draw[thick,->] ($(c.north east)+(-.4,0)$) to ($(c.north east)+(-.4,.5)$);
\draw[thick,->] ($(c.north west)+(.4,0)$) to [out=90, in= 0]
($(c.north west)+(-.1,.5)$) to [out=180, in= 90] ($(c.north west)+(-.6,0)$)
to ($(c.west)+(-.6,0)$);
\draw[thick] ($(c.west)+(-.6,0)$) to ($(c.south west)+(-.6,0)$)
to [out=-90, in= 180] ($(c.south west)+(-.1,-.5)$)
to [out=0, in= -90] ($(c.south west)+(.4,0)$);
\end{tikzpicture}
=
\begin{tikzpicture} 
[yscale=1,xscale=1,baseline=12]
\node[rectangle,draw] (c) at (1/2,1/2) {\phantom{$F_{R,T}$}};
\node[rectangle,draw] (p) at ($(c.north)+(0,1/2)$) {\phantom{$F_{R,T}$}};
\node[rectangle,draw] (i) at ($(c.south)+(0,-1/2)$) {\phantom{$F_{R,T}$}};
\node at (c) {$F_{R,T}$};
\node at (p) {$\varphi_2$};
\node at (i) {$\varphi_2^{-1}$};
\draw[thick] ($(i.south east)+(-.2,-.5)$) to ($(i.south east)+(-.2,0)$);
\draw[thick,->] ($(p.north east)+(-.2,0)$) to ($(p.north east)+(-.2,.5)$);
\draw[thick,->] ($(p.north west)+(.2,0)$) to [out=90, in= 0]
($(p.north west)+(-.3,.5)$) to [out=180, in= 90] ($(p.north west)+(-.8,0)$)
to ($(p.west)+(-.8,0)$) to ($(c.west)+(-.8,0)$);
\draw[thick] ($(c.west)+(-.8,0)$) to ($(i.west)+(-.8,0)$)
to ($(i.south west)+(-.8,0)$) to [out=-90, in= 180] ($(i.south west)+(-.3,-.5)$)
to [out=0, in= -90] ($(i.south west)+(.2,0)$);
\draw[thick] ($(i.north east)+(-.2,0)$) to ($(c.south east)+(-.2,0)$);
\draw[thick] ($(i.north west)+(.2,0)$) to ($(c.south west)+(.2,0)$);
\draw[thick] ($(c.north east)+(-.2,0)$) to ($(p.south east)+(-.2,0)$);
\draw[thick] ($(c.north west)+(.2,0)$) to ($(p.south west)+(.2,0)$);
\end{tikzpicture}
=
\begin{tikzpicture} 
[yscale=1,xscale=1,baseline=12]
\node[rectangle,draw] (c) at (1/2,1/2) {\phantom{$F_{R,T}$}};
\node at (c) {$F_{R,T}$};
\draw[thick] ($(c.south east)+(-.2,-.5)$) to ($(c.south east)+(-.2,0)$);
\draw[thick] ($(c.north east)+(-.2,0)$) to ($(c.north east)+(-.2,.5)$);
\draw[thick,->] ($(c.north west)+(.2,0)$) to [out=90, in= 0]
($(c.north west)+(-.3,.5)$) to [out=180, in= 90] ($(c.north west)+(-.8,0)$)
to ($(c.west)+(-.8,0)$);
\draw[thick] ($(c.west)+(-.8,0)$) to ($(c.south west)+(-.8,0)$)
to [out=-90, in= 180] ($(c.south west)+(-.3,-.5)$)
to [out=0, in= -90] ($(c.south west)+(.2,0)$);
\node[rectangle,draw,fill=white, anchor=south] (p)
at ($(c.east)+(-.2,.7)$) {\phantom{$\sigma^{-1}$}};
\node[rectangle,draw,fill=white, anchor=north] (i)
at ($(c.east)+(-.2,-.7)$) {\phantom{$\sigma^{-1}$}};
\node at (p) {$\sigma$};
\node at (i) {$\sigma^{-1}$};
\draw[thick] ($(i.south)+(0,-1/2)$) to ($(i.south)$);
\draw[thick,->] ($(p.north)$) to ($(p.north)+(0,1/2)$);
\end{tikzpicture}
=
\text{mRT}_{L,R}
\begin{tikzpicture}[yscale=1,xscale=1,baseline=12]
 \coordinate[] (c) at (1/2,1/2) ;
\node[rectangle,draw] (p) at ($(c)+(0,1/2)$) {\phantom{$\varphi_1^{-1}$}};
\node[rectangle,draw] (i) at ($(c)+(0,-1/2)$) {\phantom{$\varphi_1^{-1}$}};
\node at (p) {$\sigma$};
\node at (i) {$\sigma^{-1}$};
\draw[thick] ($(i.south)+(0,-1/2)$) to (i.south);\draw[thick] (i.north) to (p.south);
\draw[thick,->] (p.north) to ($(p.north)+(0,1/2)$);
\end{tikzpicture}
=
\text{mRT}_{L,R}
\begin{tikzpicture}[yscale=1,xscale=1,baseline=12]
\coordinate[] (c) at (1/2,1/2) ;
\draw[thick,] ($(c)+(0,-1)$) to (c);
\draw[thick,->] (c) to ($(c)+(0,1)$);
\end{tikzpicture}
\end{aligned}\tag*{\qed}
\end{equation*}\renewcommand{\qed}{}
\end{proof}

\begin{Remark}
 When two enhanced $R$-matrices are conjugate by $\varphi\colon V\to W$, then there is a strict ribbon equivalence $(F,G)$ between the ribbon categories~$\mathcal{V}$ and~$\mathcal{W}$ tensor generated by~$V$ and~$W$, respectively, and equipped with a braiding and ribbon structure induced by their respective enhanced $R$-matrices. The functors $F$ and $G$ are implemented by $\varphi$ and $\varphi^{-1}$ on generators.

In the case of a weak conjugacy, define the functor $F:\mathcal{V}\to\mathcal{W}$ implementing the mapping $V^{\otimes n} \to W^{\otimes n}$ by $v\mapsto \varphi_n(v)$, where $\varphi_1=\sigma$. There is a natural tensorator map:
\begin{align*} \mu_{n,m}\colon F(V^{\otimes n})\otimes F(V^{\otimes m})&\to F(V^{\otimes n+m}),
 \\
 \varphi_n(v)\otimes \varphi_m(w)&\mapsto\varphi_{n+m}(v\otimes w) .
 \end{align*}
However, $F$ is not necessarily braided monoidal with respect to $\mu$, as the following diagram may not commute
\begin{equation*}
\begin{tikzcd}
 F(V)\otimes F(V) \ar[r,"\mu_{1,1}"] \ar[d,swap,"\beta_{F(V)\otimes F(V)}"] & F(V\otimes V)
 \ar[d,"F(\beta_{V\otimes V})"]
 \\
 F(V)\otimes F(V)\ar[r,"\mu_{1,1}"]& F(V\otimes V).
 \end{tikzcd}
\end{equation*}
The composites of maps yield
\begin{gather*}
\varphi_1^{\otimes 2}(v\otimes w)
 \xmapsto{\mu_{1,1}}
 \varphi_2(v\otimes w)
 \xmapsto{F(\beta_{V\otimes V})}
 \varphi_2(R\cdot (v\otimes w)),
 \\
\varphi_1^{\otimes 2}(v\otimes w)
 \xmapsto{\beta_{F(V)\otimes F(V)}}
 R'\cdot \varphi_1^{\otimes 2}(v\otimes w)
 \xmapsto{\mu_{1,1}}
 \varphi_2\circ (\varphi_1^{-1})^{\otimes 2}\bigl(R'\cdot \varphi_1^{\otimes 2}(v\otimes w)\bigr),
\end{gather*}
which implies $\varphi_1^{\otimes 2}\cdot R = R'\cdot \varphi_1^{\otimes 2}$ rather than $\varphi_2\cdot R = R'\cdot \varphi_2$. In other words, a braided monoidal functor cannot realize a general weak conjugacy, but only one which is a standard conjugacy.
\end{Remark}

\subsection[Tensor product of R-matrices]{Tensor product of $\boldsymbol{R}$-matrices}
\label{sub.tensor}

In this subsection, we briefly review the tensor product $R' \widehat{\otimes} R''\in
\mathrm{End}\bigl((V\otimes W)^{\otimes 2}\bigr)$ of two $R$-matrices $R \in \mathrm{End}\bigl(V^{\otimes 2}\bigr)$
and $R' \in \mathrm{End}\bigl(W^{\otimes 2}\bigr)$, defined by
\begin{equation*}
R \widehat{\otimes} R' = \tau (R \otimes R') \tau,
\end{equation*}
where $\tau$ is the tensor swap $V\otimes W\to W\otimes V$ applied to the inner
factors. In other words, for all $v_i,v_j \in V$, $w_k, w_\ell \in W$, we have
\begin{equation*}
R \widehat{\otimes} R'((v_i \otimes w_k) \otimes (v_j \otimes w_\ell))=
\tau R(v_i \otimes v_j) R'(w_k \otimes w_\ell) .
\end{equation*}

It follows by the graphical calculus of the RT functor that the link invariant
of the tensor product of two enhanced $R$-matrices is the product of the
link invariants.

\section[Extending V\_1 to a link invariant]{Extending $\boldsymbol{V_1}$ to a link invariant}
\label{sec:V1}

In this section, we prove that the rigid $R$-matrix of the $V_1$ polynomial
of~\cite{GK:multi}, which we give in Appendix~\ref{sec:app1} is enhanced
and satisfies the properties \eqref{item:P1} and \eqref{item:P2} and hence leads
to a~polynomial link invariant.
In subsequent work, we will prove that this extension coincides with the
Links--Gould invariant of links, verifying a conjecture of~\cite{GK:multi}.

Our first task is to find an enhancement $R_{r}$ of the $R$-matrix
of Appendix~\ref{sec:app1}. The next lemma determines all diagonal enhancements
-- we will note later that this is a reasonable assumption.

\begin{Lemma}
\label{lem.V10}
Let $h=\mathrm{diag}(a,b,c,d) \in \mathrm{End}(W)$. Then $(R_{r},h)$ is an enhanced $R$-matrix
if and only if $r = \pm 1$ and $(a,b,c,d) = (-1, 1, 1, -1)$.
\end{Lemma}

\begin{proof}
We just have to check equations~\eqref{eq1}. Equation \eqref{eq1d}
is satisfied for any $r$ because $R_r$ is a rigid $R$-matrix~\cite{GK:multi}.
Equation~\eqref{eq1b} is equivalent to a system of four polynomial equations
with unique solution
$(a,b,c,d) = (-1, 1, 1, -1)$.
With this choice, one can check then that equation~\eqref{eq1a} is automatically
satisfied. Then, a~computation shows that the left side of equation~\eqref{eq1c}
\begin{gather*}
\bigl(R_r^{-1}\bigr)^{\circlearrowleft} \circ \bigl((\text{id}_V \otimes h) \circ R_r
\circ \bigl(h^{-1} \otimes \text{id}_V\bigr)\bigr)^{\circlearrowright}\\
\qquad{} = \text{id}_V \otimes
\text{id}_{V^{\ast}} + \bigl(1 - r^2\bigr)(1 - t_0) (e_2 \otimes e_1^\ast)^\ast
(e_4 \otimes e_3^\ast).
\end{gather*}
So $r^2 =1$. The result follows.
\end{proof}

Our next task is to show that the corresponding tangle invariant satisfies
the properties~\eqref{item:P1} and~\eqref{item:P2}.

\begin{Lemma}\label{lem.V11}
The enhanced $R$-matrix $R_{\pm 1}$ satisfies~\eqref{item:P1} and~\eqref{item:P2}.
\end{Lemma}

\begin{proof}
The proof closely follows the one of Lemma~\ref{lem.alex1}, except that the
underlying vector space of the $R$-matrix is 4-dimensional instead of
2-dimensional.

Recall that $R_{\pm1} \in \mathrm{End}(W)$ where $W$ is a 4-dimensional vector space with
$(e_1,e_2,e_3,e_4)$ a~basis for $W$. We will use not one, but two
gradings on the basis of $W$ given by
\begin{alignat*}{5}
& \mathrm{deg}_{e_2}(e_1) = 0, \qquad&& \mathrm{deg}_{e_2}(e_2) = 1, \qquad&&
\mathrm{deg}_{e_2}(e_3) = 0, \qquad&& \mathrm{deg}_{e_2}(e_4) = 1,& \\
& \mathrm{deg}_{e_3}(e_1) = 0, \qquad&& \mathrm{deg}_{e_3}(e_2) = 0, \qquad&&
\mathrm{deg}_{e_3}(e_3) = 1, \qquad&& \mathrm{deg}_{e_3}(e_4) = 1 .&
\end{alignat*}
It follows by an explicit inspection of the entries of the $R$-matrix
that the maps $(R_{\pm 1})^{\pm 1}$ are degree-preserving for both these degrees.
Actually, these gradings are not accidental. They follow from the fact that the
rank~2 Nichols algebra, which is responsible for the $R$-matrix, has natural
gradings. Moreover, a enhancement which preserves these gradings must be diagonal.
We denote the corresponding RT-functor by $F_{V_1}$ since the corresponding knot
polynomial of~\cite{GK:multi} was denoted by $V_1$.

Now we prove that $R_{\pm 1}$ satisfies \eqref{item:P1}. Fix a $(1,1)$-tangle $T$.
Using the degree-preservation of $F_{V_1,T} \in \mathrm{End}(W)$ under both gradings,
it follows that
\begin{equation}
\label{JT1}
F_{V_1,T} = \begin{pmatrix}
\alpha & 0 & 0 & 0 \\
0 & \beta & 0 & 0 \\
0 & 0 & \gamma & 0 \\
0 & 0 & 0 & \delta
\end{pmatrix} .
\end{equation}

The invariance of $F_{V_1}$ under the isotopy of
Figure~\ref{fig:egalite} implies that
\begin{equation*}
(\mathrm{id}_W \otimes F_{V_1,T}) \circ R_{\pm 1} =
R_{\pm 1} \circ (F_{V_1,T} \otimes \mathrm{id}_W) .
\end{equation*}
Inserting~\eqref{JT1} in the above identity, it follows that
$\alpha = \beta = \gamma = \delta$. Hence
$F_{V_1,T} = \alpha\,\mathrm{id}_W$, concluding the proof of the \eqref{item:P1}
identity.

Finally, we prove that $R_{1}$ satisfies \eqref{item:P2}, the proof for $R_{-1}$
involves similar computations. Fix a $(2,2)$-tangle $T$ and consider its
left and right closures $T_1$ and $T_2$, as in
Figure~\ref{fig:T12}. Using the degree-preservation of the $R$-matrix, it
follows that the operator invariant representing $T$
has the following form when written in the standard basis for $W \otimes W$:
\[
F_{V_1,T}=\left(\begin{array}{@{}cccc|cccc|cccc|cccc@{}}
a & 0 & 0 & 0 & 0 & 0 & 0 & 0 & 0 & 0 & 0 & 0 & 0 & 0 & 0 & 0 \\
0 & b & 0 & 0 & j & 0 & 0 & 0 & 0 & 0 & 0 & 0 & 0 & 0 & 0 & 0 \\
0 & 0 & d & 0 & 0 & 0 & 0 & 0 & s & 0 & 0 & 0 & 0 & 0 & 0 & 0 \\
0 & 0 & 0 & f & 0 & 0 & m & 0 & 0 & u & 0 & 0 & B & 0 & 0 & 0 \\\hline
0 & c & 0 & 0 & k & 0 & 0 & 0 & 0 & 0 & 0 & 0 & 0 & 0 & 0 & 0 \\
0 & 0 & 0 & 0 & 0 & l & 0 & 0 & 0 & 0 & 0 & 0 & 0 & 0 & 0 & 0 \\
0 & 0 & 0 & g & 0 & 0 & n & 0 & 0 & v & 0 & 0 & C & 0 & 0 & 0 \\
0 & 0 & 0 & 0 & 0 & 0 & 0 & q & 0 & 0 & 0 & 0 & 0 & F & 0 & 0 \\\hline
0 & 0 & e & 0 & 0 & 0 & 0 & 0 & t & 0 & 0 & 0 & 0 & 0 & 0 & 0 \\
0 & 0 & 0 & h & 0 & 0 & o & 0 & 0 & w & 0 & 0 & D & 0 & 0 & 0 \\
0 & 0 & 0 & 0 & 0 & 0 & 0 & 0 & 0 & 0 & y & 0 & 0 & 0 & 0 & 0 \\
0 & 0 & 0 & 0 & 0 & 0 & 0 & 0 & 0 & 0 & 0 & z & 0 & 0 & H & 0 \\\hline
0 & 0 & 0 & i & 0 & 0 & p & 0 & 0 & x & 0 & 0 & E & 0 & 0 & 0 \\
0 & 0 & 0 & 0 & 0 & 0 & 0 & r & 0 & 0 & 0 & 0 & 0 & G & 0 & 0 \\
0 & 0 & 0 & 0 & 0 & 0 & 0 & 0 & 0 & 0 & 0 & A & 0 & 0 & I & 0 \\
0 & 0 & 0 & 0 & 0 & 0 & 0 & 0 & 0 & 0 & 0 & 0 & 0 & 0 & 0 & J \\
 \end{array}\right),
\]
where $a, b, \dots, J$ are $36$ parameters. The tangle isotopies of
the left- and right-hand sides of Figure~\ref{fig:egalite2}
imply the equations
\begin{gather*}
 R_1^2 \circ F_{V_1,T} - F_{V_1,T} \circ R_1^2 = 0, \\
(R_1 \otimes \text{id}_W) \circ (\text{id}_W \otimes R_1) \circ
( F_{V_1,T} \otimes \text{id}_W) - (\text{id}_W \otimes F_{V_1,T})
\circ (R_1 \otimes \text{id}_W) \circ (\text{id}_W \otimes R_1) = 0 .
\end{gather*}
This is a linear system in $36$ variables $a,b,\dots,J$ with coefficients
in $\mathbb{Q}(t_0,t_1)$. It has rank $3$ and solving, we find that we can express
all variables as $\mathbb{Q}(t_0,t_1)$-linear combinations of $(E,i,p)$. We write the 33 other variables as inner products of $(E,i,p)$ with some triple:
\begin{alignat*}{3}
& a = \left( 1, 1, \frac{-1+t_1}{t_1} \right) \cdot (E,i,p) ,\qquad &&
s = \left(0, \frac{1}{t_1}, \frac{1}{t_1} \right) \cdot (E,i,p) ,&\\
&b = \left(1, \frac{-1+t_0}{t_0}, \frac{(-1+t_1)(-2+t_0)}{t_1(-1+t_0)}\right) \cdot (E,i,p) ,\quad &&
t = \left(1, 0, -\frac{1}{t_1}\right) \cdot (E,i,p) ,& \\
&c = \left(0, 1, \frac{(-1+t_1)t_0}{t_1(-1+t_0)}\right) \cdot (E,i,p) ,\qquad &&
u = \left(0, \frac{-1+t_0}{t_1 t_0}, \frac{1}{t_1}\right) \cdot (E,i,p) ,&\\
&d = \left(1, \frac{-1+t_1}{t_1}, \frac{-2+t_1}{t_1}\right) \cdot (E,i,p) ,\qquad &&
v = \left(0, \frac{1}{t_1}, \frac{-1+t_1 t_0}{t_1^2(-1+t_0)}\right) \cdot (E,i,p) ,&\\
&e = ( 0, 1, 1 ) \cdot (E,i,p) ,\qquad &&
w = \left(1, 0, \frac{-2+t_1+t_0}{t_1(-1+t_0)} \right) \cdot (E,i,p) ,&\\
&f = \left(1, 2 - \frac{1}{t_1} - \frac{1}{t_0}, 2 - \frac{2}{t_1} \right) \cdot (E,i,p) ,\qquad &&
x = \left(0, 0, -\frac{1}{t_1}\right) \cdot (E,i,p) ,&\\
&g = \left(0, \frac{-1+t_1}{t_1}, \frac{(-1+t_1)t_0}{t_1(-1+t_0)}\right) \cdot (E,i,p) ,\qquad &&
y = \left(1, -\frac{1}{t_1}, -\frac{2}{t_1}\right) \cdot (E,i,p) ,&\\
&h = \left(0, 1-t_0, -\frac{(-1+t_1)t_0}{-1+t_0}\right) \cdot (E,i,p) ,\qquad &&
z = \left(1, \frac{-1+t_1}{t_1}, \frac{-2+t_1}{t_1} \right) \cdot (E,i,p) ,& \\
&j = \left(0, \frac{1}{t_0}, \frac{-1+t_1}{t_1(-1+t_0)}\right) \cdot (E,i,p) ,\qquad &&
A = \left(0, -\frac{1}{t_1}, -\frac{1}{t_1}\right) \cdot (E,i,p) ,&\\
&k = \left(1, 0, \frac{1-t_1}{t_1(-1+t_0)} \right) \cdot (E,i,p) ,\qquad &&
B = \left(0, \frac{1}{t_1 t_0}, 0\right) \cdot (E,i,p) ,&\\
&l = \left(1, -\frac{1}{t_0}, -\frac{2(-1+t_1)}{t_1(-1+t_0)} \right) \cdot (E,i,p) ,\qquad &&
C = \left(0, 0, \frac{1-t_1}{t_1^2(-1+t_0)}\right) \cdot (E,i,p) ,& \\
&m = \left(0, -1 + \frac{1}{t_0}, -1 \right) \cdot (E,i,p) ,\qquad &&
D = \left(0, 0, \frac{-1+t_1}{t_1(-1+t_0)}\right) \cdot (E,i,p) ,&\\
&n = \left(1, 0, -\frac{-2+t_1+t_0}{t_1(-1+t_0)} \right) \cdot (E,i,p) ,\qquad &&
F = \left(0, -\frac{1}{t_1 t_0}, \frac{1-t_1}{t_1^2(-1+t_0)} \right) \cdot (E,i,p) ,&\\
&o = \left(0, t_1, \frac{-1+t_1 t_0}{-1+t_0} \right) \cdot (E,i,p) ,\qquad&&
G = \left(1, 0, \frac{1-t_1}{t_1(-1+t_0)}\right) \cdot (E,i,p) ,&\\
&q = \left(1, \frac{-1+t_0}{t_0}, \frac{(-1+t_1)(-2+t_0)}{t_1(-1+t_0)}\right) \cdot (E,i,p) ,\quad &&
H = (0, -1, -1) \cdot (E,i,p) ,&\\
&r = \left(0, -t_1, \frac{t_0-t_1 t_0}{-1+t_0} \right) \cdot (E,i,p) ,\qquad &&
I = \left(1, 0, -\frac{1}{t_1}\right) \cdot (E,i,p) ,&\\
&&& J = \left(1, 1, \frac{-1+t_1}{t_1}\right) \cdot (E,i,p) .&
\end{alignat*}
Using this and the definition of $\mathrm{tr}_2$ from equation~\eqref{tr2}, we can now
compute $F_{V_1,T_1}$ and $F_{V_1,T_2}$, and we find
\begin{gather*}
F_{V_1,T_1}  = \mathrm{tr}_{2} ((\mathrm{id}_W \otimes h) \circ F_{V_1,T})
= \mathrm{diag}(x_1,x_2,x_3,x_4), \\
F_{V_1,T_2}  = \mathrm{tr}_{1} (F_{V_1,T}
\circ (h^{-1} \otimes \mathrm{id}_W)) = \mathrm{diag}(x_4,x_2,x_3,x_1),
\end{gather*}
where
\begin{gather*}
 x_1  = - i - 2 p + \frac{p}{t_1} - \frac{2 p}{t_0-1}
+ \frac{2 p}{(t_0-1)t_1} + \frac{p t_0}{t_0-1} - \frac{p t_0}{(t_0-1)t_1},
\\
x_2  = - i + \frac{p}{(t_0-1)t_1} - \frac{p t_0}{t_0-1},
\\
x_3  = - i - p + \frac{p}{t_1} - \frac{p}{t_0-1} + \frac{2 p}{(t_0-1)t_1}
- \frac{p t_0}{(t_0-1)t_1},
\\
x_4  = -E - i - p - \frac{E}{t_0-1} - \frac{p}{t_0-1} + \frac{p}{(t_0-1)t_1}
+ \frac{E t_0}{t_0-1}.
\end{gather*}
This implies that $F_{V_1,T_1}=F_{V_1,T_2}$ and this
concludes the proof of \eqref{item:P2}.
\end{proof}

This concludes the proof of Theorem~\ref{thm.1} for $V_1$.

\section[Extending and identifying Lambda\_1 to a link invariant]{Extending and identifying $\boldsymbol{\Lambda_1}$ to a link invariant}
\label{sec:L1}

In this short section, we extend the $\Lambda_1$ polynomial of knots to links, and
what is more, we identify it with a known polynomial invariant.
We achieve both tasks at once by proving that the $\Lambda_1$ $R$-matrix is
simply conjugate to the product of two Alexander $R$-matrices.

To begin with, the $R$-matrix of $\Lambda_1$ given in Appendix~\ref{sec:app2},
has the following enhancement:
\begin{gather*}
\Tilde{R}_{\Lambda_1} = \frac{1}{t_0^{1/2} t_1^{1/2}} R_{\Lambda_1}
\in \mathrm{End}(X \otimes X),
\qquad h_{\Lambda_1} = t_0^{1/2} t_1^{1/2} \begin{pmatrix}
1 & 0 & 0 & 0 \\
0 & -1 & 0 & 0 \\
0 & 0 & -1 & 0 \\
0 & 0 & 0 & 1
\end{pmatrix} \in \mathrm{End}(X) .
\end{gather*}
The proof that this is indeed an enhancement follows by an explicit calculation.

We now discuss a conjugation of the above $R$-matrix with the product of two
Alexander matrices, each as in equation~\eqref{alexRh}.

Consider a 2-dimensional vector space $Y$ with basis $(y_0 , y_1)$. Let us define
the following isomorphisms:
\begin{alignat*}{3}
&\theta\colon \
 Y^{\otimes 2} \to X, \qquad&&
 \theta(y_0 \otimes y_0) = x_1,\qquad
 \theta(y_0 \otimes y_1) = x_2,& \\
 &&& \theta(y_1 \otimes y_0) = x_3,\qquad
 \theta(y_1 \otimes y_1) = x_4,&\\
&\tau\colon \
 Y^{\otimes 4} \to Y^{\otimes 4}, \qquad&&
 \tau(a \otimes b \otimes c \otimes d)= a \otimes c \otimes b \otimes d .&
\end{alignat*}

\begin{Lemma}\label{lem.RR}
We have
\begin{align*}
h_{\Lambda_1}  = \theta \circ (h_{t_1} \otimes h_{t_0}) \circ \theta^{-1},\qquad
\Tilde{R}_{\Lambda_1}   = (\theta \otimes \theta) \circ \tau
\circ (R_{t_1} \otimes R_{t_0}) \circ \tau \circ \bigl(\theta^{-1} \otimes \theta^{-1}\bigr) .
\end{align*}
\end{Lemma}

\begin{proof}
This is a direct matrix calculation.
\end{proof}

The discussion of Section~\ref{sub.conj} and the above lemma
implies that the $R$-matrices of $\Lambda_1$ and of the square of the Alexander
polynomial are conjugate, hence the two link invariants are equal.

\section[Extending and identifying Lambda\_\{-1\} to a link invariant]{Extending and identifying $\boldsymbol{\Lambda_{-1}}$ to a link invariant}
\label{sec:Lm1}

\subsection[Extension of Lambda\_\{-1\} to links]{Extension of $\boldsymbol{\Lambda_{-1}}$ to links}

We extend the $\Lambda_{-1}$ polynomial of knots to links, and identify it with the
operator invariant determined by the quantum group of $\mathfrak{sl}_3$ at a fourth root
of unity studied in~\cite{Harper}. Thus
concluding the proof of Theorems~\ref{thm.1} and~\ref{thm.2}.
Unlike the $R$-matrix $\Lambda_1$ and product of Alexander polynomials, the
$R$-matrices of $\Lambda_{-1}$ and $\Delta_{\mathfrak{sl}_3}$ are not conjugate. Yet, we
will show that they are weakly-conjugate.

Let $r^2=1$ and $\zeta=r\sqrt{-1}$ be a primitive fourth root of unity. We consider
an 8-dimensional vector space $Z$ with 2-variable Laurent polynomial coefficients
over $\mathbb{C}$. The $R$-matrix $R_{\Lambda_{-1}}$ given in Appendix \ref{sec:app3} is
expressed in the standard tensor product basis $z_i\otimes z_j$ for $1\leq i,j\leq 8$.
Moreover, $R_{\Lambda_{-1}}$ determines an enhanced $R$-matrix
$(\tilde{R}_{\Lambda_{-1}},h_{\Lambda_{-1}})$, where
\begin{gather*}
\tilde{R}_{\Lambda_{-1}}=\frac{1}{st}R_{\Lambda_{-1}}\in\mathrm{End}\bigl(Z^{\otimes 2}\bigr),\qquad
h_{\Lambda_{-1}}=st \cdot\mathrm{diag}(1,-1,-1,1,1,-1,-1,1) \in \mathrm{End}(Z) .
\end{gather*}

\subsection[The Delta\_\{sl\_3\} link invariant]{The $\boldsymbol{\Delta_{\mathfrak{sl}_3}}$ link invariant}

The $\Delta_{\mathfrak{sl}_3}$ link invariant comes from an enhanced $R$-matrix
$(R_{\mathfrak{sl}_3},h_{\mathfrak{sl}_3})$ which is given in Appendix~\ref{sec:app4} as determined
by a certain quantum group representation $V$.
This enhanced $R$-matrix satisfies~\eqref{item:P1} because of Schur's lemma, i.e.,
$V$ is a simple module. This representation is also known to be ambidextrous in the
sense of~\cite{GPT}, as shown in~\cite{Harper} following \cite{GP13}, and therefore
satisfies \eqref{item:P2}. Hence it leads to the link invariant $\Delta_{\mathfrak{sl}_3}$
of~\cite{Harper}.

To each vector in $V$, we assign a weight in the $\mathsf{A}_2$ weight lattice
\begin{alignat*}{5}
&\mathrm{wt}(v_1)=0, \qquad\! && \mathrm{wt}(v_2)=\alpha_1, \qquad\! &&
\mathrm{wt}(v_3)=\alpha_2, \qquad\! && \mathrm{wt}(v_4)=\alpha_1+\alpha_2,&
\\
&\mathrm{wt}(v_5)=\alpha_1+\alpha_2, \qquad\! &&
\mathrm{wt}(v_6)=2\alpha_1+\alpha_2, \qquad\! &&
\mathrm{wt}(v_7)=\alpha_1+2\alpha_2, \qquad\! && \mathrm{wt}(v_8)=2\alpha_1+2\alpha_2.&
\end{alignat*}
The associated degree $\mathrm{deg}_{i}$ of the basis vector $v_j$ is the coefficient
of~$\alpha_i$ in $\mathrm{wt}(v_j)$. Extend the weight function so that it is additive over
tensor products, that is $\mathrm{wt}(v_i\otimes v_j)=\mathrm{wt}(v_i)+\mathrm{wt}(v_j)$.

The $R$-matrix respects this weight grading on $V\otimes V$. Therefore, any given
weight determines an invariant subspace. This justifies our presentation of the
$R$-matrices in Appendices~\ref{sec:app3} and~\ref{sec:app4}. Moreover,
$R_{\mathfrak{sl}_3}$ and $\tilde{R}_{\Lambda_{-1}}$ have the same support.

\subsection[Weak conjugacy between Delta\_\{sl\_3\} and Lambda\_\{-1\}]{Weak conjugacy between $\boldsymbol{\Delta_{\mathfrak{sl}_3}}$ and $\boldsymbol{\Lambda_{-1}}$}

In this subsection, we show that the enhanced $R$-matrices $R_{\mathfrak{sl}_3}$
and $\tilde{R}_{\Lambda_{-1}}$ are weakly-conjugate, hence the corresponding
link invariants coincide. To do so, we define the following invertible linear maps:
\begin{gather*}
\sigma \colon \ V \to Z, \qquad (z_1 ,z_2 ,z_3 ,z_4 ,z_5 ,z_6 ,z_7 ,z_8)
\mapsto (v_1 ,v_2 ,v_3 ,v_4 ,\zeta v_5 ,\zeta v_6 ,\zeta v_7 , \zeta v_8),
\\
\nu \colon \ V\to V, \qquad v_i \mapsto t_1^{-\mathrm{deg}_{1}(v_i)}t_2^{-\mathrm{deg}_{2}(v_i)}v_i,
\\
\gamma \colon \ V \otimes V \to V \otimes V,\qquad
v_i\otimes v_j \mapsto \zeta^{\mathrm{deg}_{2}(v_i)\mathrm{deg}_{1}(v_j)}v_i\otimes v_j .
\end{gather*}
Let $\varphi=(\sigma\otimes\sigma)\circ (\nu\otimes \mathrm{id}_V)\circ \gamma \in \mathrm{Hom}\bigl(V^{\otimes 2},Z^{\otimes 2}\bigr)$.

\begin{Lemma}\label{lem:R2}
We have
\begin{gather*}
 \varphi\circ {R}_{\mathfrak{sl}_3}\circ\varphi^{-1}=  \bigl\lbrack
\tilde{R}_{\Lambda_{-1}} \bigr\rbrack_{t=t_1^{-2}, s=t_2^{-2}} .
\end{gather*}
\end{Lemma}
The equality of the lemma can be proven by comparing the two endomorphisms on invariant subspaces
of $V\otimes V$.
Let
\begin{gather*}
 \varphi_n=\biggl(\bigotimes_{i=1}^n \sigma\circ\nu^{n-i}\biggr)\circ
 \biggl(\prod_{1\leq i<j\leq n}(\gamma)_{i,j}\biggr)
\end{gather*}
In this way, $\varphi=\varphi_2$ is the intertwiner between $R$-matrices given in
Lemma~\ref{lem:R2}. Setting $\nu_{n-1}=\bigotimes_{i=1}^{n-1} \nu$ and
$\gamma_n=\prod_{1\leq i<n} (\gamma)_{i,n}$ yields the identity
$\varphi_n=(\varphi_{n-1}\otimes \mathrm{id}_Z)\circ(\nu_{n-1}\otimes\sigma)\circ \gamma_n$
in accordance with equation \eqref{eq.phin}. Equation~\eqref{eq.nun} is satisfied
for our choice of $\nu_{n-1}$ because ${\mathrm{tr}_n((1\otimes h)\circ F_T)}$ is a
weight-preserving linear map and $\nu_{n-1}$ is constant on weight spaces of~$V^{\otimes (n-1)}$. The endomorphism $\gamma_n$ is constant on weight spaces determined
by the weights of the first $n-1$ and the last tensor factors taken as a pair.
Since~$F_T$ and $h$ are both weight preserving and the $n$-th partial trace
operation preserves weight in the $n$-th tensor factor, equation~\eqref{eq.gn} is
then satisfied. Similarly, equation~\eqref{eq.g2} holds. Clearly, \eqref{eq.sn} and \eqref{eq.n1} hold for our choices of $\sigma$ and $\nu_1$, under the
specialization of parameters. Thus, $\varphi_n$ satisfies property~\eqref{item:BC2}.

\begin{Lemma}\label{lem.Rn}
For each $n\geq 2$ and $1\leq k\leq n-1$, there is an equality of morphisms on
$\mathrm{End}(Z^{\otimes n})$
\begin{gather*}
\varphi_n\circ({R}_{\mathfrak{sl}_3})_{k,k+1}\circ\varphi_n^{-1}
= \bigl\lbrack \bigl(\tilde{R}_{\Lambda_{-1}}\bigr)_{k,k+1} \bigr\rbrack_{t=t_1^{-2}, s=t_2^{-2}} .
\end{gather*}
Therefore, $\varphi_n$ satisfies~\eqref{item:BC1} for weak conjugacy.
\end{Lemma}

\begin{proof}Any morphism which does not act on tensor factors $k$ or $k+1$ commutes with
$({R}_{\mathfrak{sl}_3})_{k,k+1}$, thereby reducing our considerations to conjugation
of $({R}_{\mathfrak{sl}_3})_{k,k+1}$ by
\begin{gather}
((\sigma \otimes \sigma)
\circ (\nu\otimes \mathrm{id}_V)\circ \gamma)_{k,k+1}\nonumber
\\
\qquad{} \circ(\nu^{n-(k+1)}\otimes \nu^{n-(k+1)})_{k,k+1}\circ
\biggl(\prod_{1\leq i<k}(\gamma)_{i,k}(\gamma)_{i,k+1}\biggr)\circ
\biggl(\prod_{k+1<j \leq n}(\gamma)_{k,j}(\gamma)_{k+1,j}\biggr).\!\!\!\label{eq:Rcomp1}
\end{gather}
The second line of morphisms in \eqref{eq:Rcomp1} commutes with
$(R_{\mathfrak{sl}_3})_{k,k+1}$. This is because the morphisms depend only on the weight of
the vectors in the $k$-th and $(k+1)$-st tensor factors and $R$-matrix preserves
these weights. Now applying Lemma~\ref{lem:R2} to the conjugacy by $(\varphi)_{k,k+1}$,
the first line of~\eqref{eq:Rcomp1}, and evaluating the change of parameters,
we obtain $\tilde{R}_{\Lambda_{-1}}$.
\end{proof}

Since the enhanced $R$-matrices $R_{\mathfrak{sl}_3}$ and $\tilde{R}_{\Lambda_{-1}}$ are
weakly-conjugate, it follows from Lemma~\ref{lem.RTequal} that the corresponding
link invariants $\Lambda_{-1}$ and $\Delta_{\mathfrak{sl}_3}$ are equal, after a change of
variables $(t,s) \to \bigl(t_1^{-2},t_2^{-2}\bigr)$.
This concludes equation~\eqref{Lam1} of Theorem~\ref{thm.2}.

\appendix

\section[Explicit R-matrices]{Explicit $\boldsymbol{R}$-matrices}

In this appendix, we give the four $R$-matrices that appear in our paper.

\subsection[The R-matrix of the V\_1-polynomial]{The $\boldsymbol{R}$-matrix of the $\boldsymbol{V_1}$-polynomial}
\label{sec:app1}

In this appendix, we give the enhanced $R$-matrix of the $V_1$-polynomial
whose explicit computation was discussed in~\cite{GK:multi}. We fix two variables
$t_0$, $t_1$ and consider a 4-dimensional vector space~$W$ with
basis $(e_1,e_2,e_3,e_4)$. We define a map $R_{r} \in \mathrm{End}{(W\otimes W)}$ that
depends on an additional free variable~$r$ and can be presented in matrix form
as follows, when written in the standard basis
$(e_1 \otimes e_1, e_1 \otimes e_2, e_1 \otimes e_3,
e_1 \otimes e_4, e_2 \otimes e_1, e_2 \otimes e_2, \dots)$:
\[
\resizebox{0.99\textwidth}{!}{$
R_{r}=\left(\begin{array}{@{}c@{\,}c@{\,}c@{\,}c@{\,}|@{\,}c@{\,}c@{\,}c@{\,}c@{\,}|@{\,}c@{\,}c@{\,}c@{\,}c@{\,}|@{\,}c@{\,}c@{\,}c@{\,}c@{}}
-1 & 0 & 0 & 0 & 0 & 0 & 0 & 0 & 0 & 0 & 0 & 0 & 0 & 0 & 0 & 0 \\
0 & 0 & 0 & 0 & -1 & 0 & 0 & 0 & 0 & 0 & 0 & 0 & 0 & 0 & 0 & 0 \\
0 & 0 & 0 & 0 & 0 & 0 & 0 & 0 & -1 & 0 & 0 & 0 & 0 & 0 & 0 & 0 \\
0 & 0 & 0 & 0 & 0 & 0 & 0 & 0 & 0 & 0 & 0 & 0 & -1 & 0 & 0 & 0 \\\hline
0 & -t_{0} & 0 & 0 & t_{0}-1 & 0 & 0 & 0 & 0 & 0 & 0 & 0 & 0 & 0 & 0 & 0 \\
0 & 0 & 0 & 0 & 0 & t_{0} & 0 & 0 & 0 & 0 & 0 & 0 & 0 & 0 & 0 & 0 \\
0 & 0 & 0 & 0 & 0 & 0 & 0 & 0 & 0 & -r^{-1}t_{1}^{-1} & 0 & 0 & t_{1}^{-1}-1
& 0 & 0 & 0 \\
0 & 0 & 0 & 0 & 0 & 0 & 0 & 0 & 0 & 0 & 0 & 0 & 0 & r^{-1}t_{1}^{-1}
& 0 & 0 \\\hline
0 & 0 & -t_{1} & 0 & 0 & 0 & 0 & 0 & t_{1}-1 & 0 & 0 & 0 & 0 & 0 & 0 & 0 \\
0 & 0 & 0 & 0 & 0 & 0 & -r t_{1} & 0 & 0 & 0 & 0 & 0 & r(t_{1}-1) & 0 & 0 & 0 \\
0 & 0 & 0 & 0 & 0 & 0 & 0 & 0 & 0 & 0 & t_{1} & 0 & 0 & 0 & 0 & 0 \\
0 & 0 & 0 & 0 & 0 & 0 & 0 & 0 & 0 & 0 & 0 & 0 & 0 & 0 & r t_{1} & 0 \\\hline
0 & 0 & 0 & -t_{0}t_{1} & 0 & 0 & (t_{0}-1)t_{1} & 0 & 0 & r^{-1}(1-t_{0})
& 0 & 0 & t_{0}+t_{1}-2 & 0 & 0 & 0 \\
0 & 0 & 0 & 0 & 0 & 0 & 0 & r t_{0}t_{1} & 0 & 0 & 0 & 0 & 0 & t_{0}-1 & 0 & 0 \\
0 & 0 & 0 & 0 & 0 & 0 & 0 & 0 & 0 & 0 & 0 & r^{-1} & 0 & 0 & t_{1}-1 & 0 \\
0 & 0 & 0 & 0 & 0 & 0 & 0 & 0 & 0 & 0 & 0 & 0 & 0 & 0 & 0 & -1 \\
\end{array}\right)$}.
\]
The basis can be reordered, which presents the matrix in block form.
Here $(i,j)$ stands for $e_i \otimes e_j$:
\[
\resizebox{0.99\textwidth}{!}{$
R_{r}=\begin{blockarray}{c@{\,}c@{\,}c@{\,}c@{\,}c@{\,}c@{\,}c@{\,}c@{\,}c@{\,}c@{\,}c@{\,}c@{\,}c@{\,}c@{\,}c@{\,}c@{\,}c@{\,}c@{\,}c@{}}
& \matindex{(1,1)} & \matindex{(1,2)} & \matindex{(2,1)} & \matindex{(1,3)}
&\matindex{(3,1)} & \matindex{(2,2)} &\matindex{(3,3)} & \matindex{(1,4)}
&\matindex{(2,3)} & \matindex{(3,2)} &\matindex{(4,1)} & \matindex{(2,4)}
&\matindex{(4,2)} & \matindex{(3,4)} &\matindex{(4,3)} & \matindex{(4,4)}& \\
\begin{block}{(@{\,}c@{\,}c@{\,}c@{\,}c@{\,}c@{\,}c@{\,}c@{\,}c@{\,}c@{\,}c@{\,}c@{\,}c@{\,}c@{\,}c@{\,}c@{\,}c@{\,}c@{\,}c@{\,})@{\,\,\,}c@{\,}}
& -1 & . & . & . & . & . & . & . & . & . & . & . & . & . & . & . &
&\matindex{(1,1)} \\
& . & 0 & -1 & . & . & . & . & . & . & . & . & . & . & . & . & . &
&\matindex{(1,2)} \\
& . & -t_{0} & t_{0}-1 & . & . & . & . & . & . & . & . & . & . & . & . & . &
&\matindex{(2,1)} \\
& . & . & . & 0 & -1 & . & . & . & . & . & . & . & . & . & . & . &
&\matindex{(1,3)} \\
& . & . & . & -t_{1} & t_{1}-1 & . & . & . & . & . & . & . & . & . & . & . &
&\matindex{(3,1)} \\
& . & . & . & . & . & t_{0} & . & . & . & . & . & . & . & . & . & . &
&\matindex{(2,2)} \\
& . & . & . & . & . & . & t_{1} & . & . & . & . & . & . & . & . & . &
&\matindex{(3,3)} \\
& . & . & . & . & . & . & . & 0 & 0 & 0 & -1 & . & . & . & . & . &
&\matindex{(1,4)} \\
& . & . & . & . & . & . & . & 0 & 0 & -r^{-1} t_{1}^{-1} & t_{1}^{-1}-1
& . & . & . & . & . & &\matindex{(2,3)} \\
& . & . & . & . & . & . & . & 0 & -r t_{1} & 0 & r(t_{1}-1) & . & . & . & .
& . & &\matindex{(3,2)} \\
& . & . & . & . & . & . & . & -t_{0}t_{1} & t_{1}(t_{0}-1) & (1-t_{0})r^{-1}
&t_{0} + t_{1} -2 & . & . & . & . & . & &\matindex{(4,1)} \\
& . & . & . & . & . & . & . & . & . & . & . & 0 & r^{-1} t_{1}^{-1} & . & .
& . & &\matindex{(2,4)} \\
& . & . & . & . & . & . & . & . & . & . & . & r t_{0}t_{1} & t_{0}-1 & . & .
& . & &\matindex{(4,2)} \\
& . & . & . & . & . & . & . & . & . & . & . & . & . & 0 & r t_{1} & . &
&\matindex{(3,4)} \\
& . & . & . & . & . & . & . & . & . & . & . & . & . & r^{-1} & t_{1}-1 & .
& &\matindex{(4,3)} \\
& . & . & . & . & . & . & . & . & . & . & . & . & . & . & . & -1 &
&\matindex{(4,4)} \\
\end{block}
\end{blockarray}
$.}
\]

\subsection[The R-matrix of the Lambda\_1-polynomial]{The $\boldsymbol{R}$-matrix of the $\boldsymbol{\Lambda_{1}}$-polynomial}
\label{sec:app2}

Fix two variables $t_0$, $t_1$ and a 4-dimensional vector space $X$ over a field
of characteristic zero with basis $(x_1, x_2, x_3, x_4)$. We can write the
rigid $R$-matrix $R_{\Lambda_1}\in \mathrm{End}(X \otimes X)$ that defines polynomial
$\Lambda_1$ using the Garoufalidis--Kashaev construction in the natural basis
for $X \otimes X$, that is $(x_1 \otimes x_1, x_1 \otimes x_2, x_1 \otimes x_3,
x_1 \otimes x_4, x_2 \otimes x_1, \dots)$:
\[
\resizebox{0.99\textwidth}{!}{
$
R_{\Lambda_1}=\left(\begin{array}{@{}c@{\,}c@{\,}c@{\,}c@{\,}|@{\,}c@{\,}c@{\,}c@{\,}c@{\,}|@{\,}c@{\,}c@{\,}c@{\,}c@{\,}|@{\,}c@{\,}c@{\,}c@{\,}c@{}}
1 & 0 & 0 & 0 & 0 & 0 & 0 & 0 & 0 & 0 & 0 & 0 & 0 & 0 & 0 & 0 \\
0 & 0 & 0 & 0 & 1 & 0 & 0 & 0 & 0 & 0 & 0 & 0 & 0 & 0 & 0 & 0 \\
0 & 0 & 0 & 0 & 0 & 0 & 0 & 0 & 1 & 0 & 0 & 0 & 0 & 0 & 0 & 0 \\
0 & 0 & 0 & 0 & 0 & 0 & 0 & 0 & 0 & 0 & 0 & 0 & 1 & 0 & 0 & 0 \\\hline
0 & t_0 & 0 & 0 & 1-t_0 & 0 & 0 & 0 & 0 & 0 & 0 & 0 & 0 & 0 & 0 & 0 \\
0 & 0 & 0 & 0 & 0 & -t_0 & 0 & 0 & 0 & 0 & 0 & 0 & 0 & 0 & 0 & 0 \\
0 & 0 & 0 & 0 & 0 & 0 & 0 & 0 & 0 & t_0 r & 0 & 0 & 1-t_0 & 0 & 0 & 0 \\
0 & 0 & 0 & 0 & 0 & 0 & 0 & 0 & 0 & 0 & 0 & 0 & 0 & - t_0 r & 0 & 0 \\\hline
0 & 0 & t_1 & 0 & 0 & 0 & 0 & 0 & 1-t_1 & 0 & 0 & 0 & 0 & 0 & 0 & 0 \\
0 & 0 & 0 & 0 & 0 & 0 & t_1 r^{-1} & 0 & 0 & 0 & 0 & 0 & (1-t_1) r^{-1} & 0
& 0 & 0 \\
0 & 0 & 0 & 0 & 0 & 0 & 0 & 0 & 0 & 0 & -t_1 & 0 & 0 & 0 & 0 & 0 \\
0 & 0 & 0 & 0 & 0 & 0 & 0 & 0 & 0 & 0 & 0 & 0 & 0 & 0 & -t_1 r^{-1} & 0
\\\hline
0 & 0 & 0 & t_0 t_1 & 0 & 0 & (1-t_0) t_1 & 0 & 0 & t_0 (1-t_1) r & 0 & 0
& (1-t_0)(1-t_1) & 0 & 0 & 0 \\
0 & 0 & 0 & 0 & 0 & 0 & 0 & - t_0 t_1 r^{-1} & 0 & 0 & 0 & 0 & 0 & -t_0 (1-t_1)
& 0 & 0 \\
0 & 0 & 0 & 0 & 0 & 0 & 0 & 0 & 0 & 0 & 0 & -t_0 t_1 r & 0 & 0 & -t_1 (1-t_0)
& 0 \\
0 & 0 & 0 & 0 & 0 & 0 & 0 & 0 & 0 & 0 & 0 & 0 & 0 & 0 & 0 & t_0 t_1 \\
\end{array}\right)
$},
\]

\subsection[The R-matrix of the Lambda\_\{-1\}-polynomial]{The $\boldsymbol{R}$-matrix of the $\boldsymbol{\Lambda_{-1}}$-polynomial}
\label{sec:app3}

Fix two variables $t,s$ and an 8-dimensional vector space $Z$ with
basis $(z_1,\dots,z_8)$. The rigid $R$-matrix $R_{\Lambda_{-1}}\in \mathrm{End}(Z \otimes Z)$
that defines the invariant $\Lambda_{-1}$ can be
decomposed in terms of invariant subspaces in the natural basis for $Z \otimes Z$.
Rows and columns are labeled $(i,j)$ to designate the vector $z_i\otimes z_j$,
and $r$ is an additional variable.

Six 1-dimensional invariant subspaces:
\begin{alignat*}{4}
&\begin{blockarray}{cc}
\matindex{(1,1)}&\\
\begin{block}{(c)c}
1 & \matindex{(1,1)}
\\
\end{block}
\end{blockarray}
\qquad&&
\begin{blockarray}{cc}
\matindex{(2,2)}&\\
\begin{block}{(c)c}
-s & \matindex{(2,2)}
\\
\end{block}
\end{blockarray}
\qquad&&
\begin{blockarray}{cc}
\matindex{(3,3)}&\\
\begin{block}{(c)c}
-t & \matindex{(3,3)}
\\
\end{block}
\end{blockarray}&
\\[-2mm]
&\begin{blockarray}{cc}
\matindex{(6,6)}&\\
\begin{block}{(c)c}
-s^2 t & \matindex{(6,6)}
\\
\end{block}
\end{blockarray}
\qquad&&
\begin{blockarray}{cc}
\matindex{(7,7)}&\\
\begin{block}{(c)c}
-s t^2 & \matindex{(7,7)}
\\
\end{block}
\end{blockarray}
\qquad&&
\begin{blockarray}{cc}
\matindex{(8,8)}&\\
\begin{block}{(c)c}
s^2 t^2 & \matindex{(8,8)}
\\
\end{block}
\end{blockarray}&
\end{alignat*}
Six 2-dimensional invariant subspaces:
\begin{align*}
\begin{blockarray}{ccc}
\matindex{(1,2)}
&
\matindex{(2,1)}\\
\begin{block}{(cc)c}
0 & 1&\matindex{(1,2)}
\\
s & 1-s &\matindex{(2,1)}
\\
\end{block}
\end{blockarray}
&&
\begin{blockarray}{ccc}
\matindex{(1,3)}
&
\matindex{(3,1)}\\
\begin{block}{(cc)c}
0 & 1&\matindex{(1,3)}
\\
t & 1-t &\matindex{(3,1)}
\\
\end{block}
\end{blockarray}
 &&
\begin{blockarray}{ccc}
\matindex{(2,6)}
&
\matindex{(6,2)}\\
\begin{block}{(cc)c}
0 & r s &\matindex{(2,6)}
\\
-\frac{s^2 t}{r} & -s(1 + s t) &\matindex{(6,2)}
\\
\end{block}
\end{blockarray}
\\[-2mm]
\begin{blockarray}{ccc}
\matindex{(3,7)}
&
\matindex{(7,3)}\\
\begin{block}{(cc)c}
0 & -\frac{t}{r}&\matindex{(3,7)}
\\
r s t^2 & -t(1 + s t) &\matindex{(7,3)}
\\
\end{block}
\end{blockarray}
&&
\begin{blockarray}{ccc}
\matindex{(6,8)}
&
\matindex{(8,6)}\\
\begin{block}{(cc)c}
0 & r^2 s^2 t &\matindex{(6,8)}
\\
\frac{s^2 t^2}{r^2} & - s^2 t (1-t) &\matindex{(8,6)}
\\
\end{block}
\end{blockarray}
&&
\begin{blockarray}{ccc}
\matindex{(7,8)}
&
\matindex{(8,7)}\\
\begin{block}{(cc)c}
0 & \frac{s t^2}{r^2}&\matindex{(7,8)}
\\[1mm]
r^2 s^2 t^2 & -s t^2 (1-s) &\matindex{(8,7)}
\\
\end{block}
\end{blockarray}
\end{align*}
Six 6-dimensional invariant subspaces:
\begin{gather*}
\begin{blockarray}{ccccccc}
\matindex{(1,4)}
&
\matindex{(1,5)}
&
\matindex{(2,3)}
&
\matindex{(3,2)}
&
\matindex{(5,1)}
&
\matindex{(4,1)}\\
\begin{block}{(cccccc)c}
0&0&0&0&0&1&\matindex{(1,4)}
\\
0&0&0&0&1&0&\matindex{(1,5)}
\\
0&0&0&r s&r(1-s)&1+s&\matindex{(2,3)}
\\
0&0&-\frac{t}{r}&0&1+t&\frac{t-1}{r}&\matindex{(3,2)}
\\[1mm]
0&t s&0&s(1-t)&1-s&\frac{s(1-t)}{r}&\matindex{(5,1)}
\\[1mm]
s t&0&t(1-s)&0&r t (s-1)&1-t&\matindex{(4,1)}
\\
\end{block}
\end{blockarray}
\\
\begin{blockarray}{ccccccc}
\matindex{(1,6)}
&
\matindex{(2,4)}
&
\matindex{(2,5)}
&
\matindex{(5,2)}
&
\matindex{(4,2)}
&
\matindex{(6,1)}\\
\begin{block}{(cccccc)c}
0&0&0&0&0&1&\matindex{(1,6)}
\\
0&0&0&0&-r s&r(s-1)&\matindex{(2,4)}
\\
0&0&0&-r s&0&1+s&\matindex{(2,5)}
\\
0&0&\frac{st}{r}&-s(1+t)&\frac{s(1-t)}{r}&\frac{1+st}{r}&\matindex{(5,2)}
\\[1mm]
0&\frac{st}{r}&0&0&0&1+st&\matindex{(4,2)}
\\[1mm]
s^2 t&0&t s (1-s)&r s t (s-1)&s(1-t)&(1-s)(1+st)&\matindex{(6,1)}
\\
\end{block}
\end{blockarray}
\\
\begin{blockarray}{ccccccc}
\matindex{(1,7)}
&
\matindex{(3,5)}
&
\matindex{(3,4)}
&
\matindex{(4,3)}
&
\matindex{(5,3)}
&
\matindex{(7,1)}\\
\begin{block}{(cccccc)c}
0&0&0&0&0&1&\matindex{(1,7)}
\\
0&0&0&0&\frac{t}{r}&\frac{1-t}{r}&\matindex{(3,5)}
\\
0&0&0&\frac{t}{r}&0&1+t&\matindex{(3,4)}
\\
0&0&-r s t&-(1+s)t&rt(s-1)&-r(1+st)&\matindex{(4,3)}
\\
0&-r s t&0&0&0&1+st&\matindex{(5,3)}
\\
s t^2&0&t s (1-t)&\frac{ts(1-t)}{r}&t(1-s)&(1-t)(1+st)&\matindex{(7,1)}
\\
\end{block}
\end{blockarray}
\\
\begin{blockarray}{ccccccc}
\matindex{(2,8)}
&
\matindex{(5,6)}
&
\matindex{(4,6)}
&
\matindex{(6,4)}
&
\matindex{(6,5)}
&
\matindex{(8,2)}\\
\begin{block}{(cccccc)c}
0&0&0&0&0&r^2 s&\matindex{(2,8)}
\\
0&0&0&0&-\frac{st}{r}&\frac{s(t-1)}{r}&\matindex{(5,6)}
\\
0&0&0&-\frac{st}{r}&0&0&\matindex{(4,6)}
\\
0&0&rs^2t&0&rst(1-s)&rst(s-1)&\matindex{(6,4)}
\\
0&rs^2t&0&0&-st(1+s)&s(t-1)&\matindex{(6,5)}
\\
\frac{s^2t^2}{r^2}&rs^2t(1+t)&s^2t(t-1)&\frac{st(1+st)}{r}
&-st(1+st)&s(t-1)(1+st)&\matindex{(8,2)}
\\
\end{block}
\end{blockarray}
\\
\begin{blockarray}{ccccccc}
\matindex{(3,8)}
&
\matindex{(4,7)}
&
\matindex{(5,7)}
&
\matindex{(7,5)}
&
\matindex{(7,4)}
&
\matindex{(8,3)}\\
\begin{block}{(cccccc)c@{}}
0&0&0&0&0&\frac{t}{r^2}&\matindex{(3,8)}
\\[1mm]
0&0&0&0&rst&\frac{t(s-1)}{r}&\matindex{(4,7)}
\\[1mm]
0&0&0&rst&0&0&\matindex{(5,7)}
\\
0&0&-\frac{st^2}{r}&0&\frac{st(t-1)}{r}&\frac{st(t-1)}{r^3}&\matindex{(7,5)}
\\[1mm]
0&-\frac{st^2}{r}&0&0&-st(1+t)&-\frac{t(s-1)}{r^2}&\matindex{(7,4)}
\\[1mm]
r^2s^2t^2& rst^2(1+s) &r^2st^2(1-s)&r^3st(1+st)&r^2st(1+st)
&t(s-1)(1+st)&\matindex{(8,3)}
\\
\end{block}
\end{blockarray}
\\
\begin{blockarray}{ccccccc}
\matindex{(4,8)}
&
\matindex{(5,8)}
&
\matindex{(6,7)}
&
\matindex{(7,6)}
&
\matindex{(8,5)}
&
\matindex{(8,4)}\\
\begin{block}{(cccccc)c}
0&0&0&0&0&st&\matindex{(4,8)}
\\
0&0&0&0&st&0&\matindex{(5,8)}
\\
0&0&0&-r^3s^2t&rst(1-s)&0&\matindex{(6,7)}
\\
0&0&\frac{st^2}{r^3}&0&0&\frac{st(1-t)}{r^3}&\matindex{(7,6)}
\\[1mm]
0&s^2 t^2&-\frac{s(1+s)t^2}{r}&r^2 s^2 t(t-1)&s^2t(t-1)
&\frac{st(t-1)}{r}&\matindex{(8,5)}
\\[1mm]
s^2t^2&0&st^2(1-s)&-r^3s^2t(1+t)&rst(1-s)&s(s-1)t^2&\matindex{(8,4)}
\\
\end{block}
\end{blockarray}
\end{gather*}
\newpage

\noindent
One 10-dimensional invariant subspace:
\[
\resizebox{0.99\textwidth}{!}{
$\begin{blockarray}{ccccccccccc}
\matindex{(1,8)}
&
\matindex{(2,7)}
&
\matindex{(3,6)}
&
\matindex{(4,5)}
&
\matindex{(4,4)}
&
\matindex{(5,5)}
&
\matindex{(5,4)}
&
\matindex{(6,3)}
&
\matindex{(7,2)}
&
\matindex{(8,1)}\\[2ex]
\begin{block}{(cccccccccc)c}
0&0&0&0&0&0&0&0&0&1&\matindex{(1,8)}
\\
0&0&0&0&0&0&0&0&-r^2s&1-s&\matindex{(2,7)}
\\
0&0&0&0&0&0&0&-\frac{t}{r^2}&0&\frac{t-1}{r^2}&\matindex{(3,6)}
\\[1mm]
0&0&0&0&0&0&-st&\frac{t(s+1)}{r}&0&\frac{1-t}{r}&\matindex{(4,5)}
\\
0&0&0&0&-st&0&0&t(s-1)&0&(1-s)t&\matindex{(4,4)}
\\
0&0&0&0&0&-st&0&0&s(t-1)&\frac{s(t-1)}{r^2}&\matindex{(5,5)}
\\[1mm]
0&0&0&-st&0&0&0&0&-rs(1+t)&\frac{1-s}{r}&\matindex{(5,4)}
\\[1mm]
0&0&-r^2s^2t&-rst(s+1)&0&r^2st(s-1)&0&0&-r^2s(1+st)&
(1-s)(1+st)&\matindex{(6,3)}
\\
0&-\frac{st^2}{r^2}&0&0&\frac{st(t-1)}{r^2}&0&\frac{st(t+1)}{r}
&-\frac{t(1+st)}{r^2}&0& \frac{(t-1)(1+st)}{r^2}&\matindex{(7,2)}
\\[1mm]
s^2t^2&st^2(1-s)&r^2s^2t(t-1)&
rs^2t(t-1)&
st(1-t)
&
r^2st(s-1)&
rst^2(s-1)&
\begin{matrix}
t(1-s)\\
\times (1+st)
\end{matrix}&
\begin{matrix}
r^2s(t-1)\\
\times (1+st)
\end{matrix}&
\begin{matrix}
(s-1)(t-1)\\
\times
(1+st)
\end{matrix}&\matindex{(8,1)}
\\
\end{block}
\end{blockarray}
$
}
\]

\subsection[The R-matrix of Delta\_\{sl\_3\}]{The $\boldsymbol{R}$-matrix of $\boldsymbol{\Delta_{\mathfrak{sl}_3}}$}
\label{sec:app4}

Fix two variables $t_1$, $t_2$ and $\zeta=\sqrt{-1}$. Consider an 8-dimensional
vector space $V=V(t_1,t_2)$ with basis $(v_1,\dots,v_8)$.

Next define $R_{\mathfrak{sl}_3}\in\mathrm{End}{(V\otimes V)}$ in the standard tensor product
basis $(v_i\otimes v_j)$. Given that $V\otimes V$ is a 64-dimensional vector space,
to ease the presentation of the enhanced $R$-matrix $(R_{\mathfrak{sl}_3},h_{\mathfrak{sl}_3})$
we define its action in terms of blocks
(invariant subspaces). Once again, rows and columns are labeled by pairs $(i,j)$
to indicate the basis vector $v_i\otimes v_j$.

We give the enhancement:
\[
h_{\mathfrak{sl}_3} =t_1^{-2}t_2^{-2}\cdot \mathrm{diag}(1,-1,-1,1,1,-1,-1,1) .
\]
Six 1-dimensional invariant subspaces:
\begin{gather*}
\begin{blockarray}{cc}
\matindex{(1,1)}&\\
\begin{block}{(c)c}
t_1^2t_2^2 & \matindex{(1,1)}
\\
\end{block}
\end{blockarray}
\qquad
\begin{blockarray}{cc}
\matindex{(2,2)}&\\
\begin{block}{(c)c}
-t_2^2 & \matindex{(2,2)}
\\
\end{block}
\end{blockarray}
\qquad
\begin{blockarray}{cc}
\matindex{(3,3)}&\\
\begin{block}{(c)c}
-t_1^2 & \matindex{(3,3)}
\\
\end{block}
\end{blockarray}
\\
\begin{blockarray}{cc}
\matindex{(6,6)}&\\
\begin{block}{(c)c}
-t_1^{-2} & \matindex{(6,6)}
\\
\end{block}
\end{blockarray}
\qquad
\begin{blockarray}{cc}
\matindex{(7,7)}&\\
\begin{block}{(c)c}
-t_2^{-2} & \matindex{(7,7)}
\\
\end{block}
\end{blockarray}
\qquad
\begin{blockarray}{cc}
\matindex{(8,8)}&\\
\begin{block}{(c)c}
t_1^{-2}t_2^{-2} & \matindex{(8,8)}
\\
\end{block}
\end{blockarray}
\end{gather*}
Six 2-dimensional invariant subspaces:
\begin{gather*}
t_2^2
\begin{blockarray}{ccc}
\matindex{(1,2)}
&
\matindex{(2,1)}\\
\begin{block}{(cc)c}
0 & t_1&\matindex{(1,2)}
\\
t_1 & t_1^2-1 &\matindex{(2,1)}
\\
\end{block}
\end{blockarray}
\qquad
t_1^2
\begin{blockarray}{ccc}
\matindex{(1,3)}
&
\matindex{(3,1)}\\
\begin{block}{(cc)c}
0 & t_2&\matindex{(1,3)}
\\
t_2 & t_2^2-1 &\matindex{(3,1)}
\\
\end{block}
\end{blockarray}
\\
-t_1^{-1}t_2
\begin{blockarray}{ccc}
\matindex{(2,6)}
&
\matindex{(6,2)}\\
\begin{block}{(cc)c}
0 & \zeta&\matindex{(2,6)}
\\
\zeta & t_1t_2+t_1^{-1}t_2^{-1} &\matindex{(6,2)}
\\
\end{block}
\end{blockarray}
\qquad
-t_1t_2^{-1}
\begin{blockarray}{ccc}
\matindex{(3,7)}
&
\matindex{(7,3)}\\
\begin{block}{(cc)c}
0 & \zeta&\matindex{(3,7)}
\\
\zeta & t_1t_2+t_1^{-1}t_2^{-1} &\matindex{(7,3)}
\\
\end{block}
\end{blockarray}
\\
 -t_1^{-2}
\begin{blockarray}{ccc}
\matindex{(6,8)}
&
\matindex{(8,6)}\\
\begin{block}{(cc)c}
0 & t_2^{-1}&\matindex{(6,8)}
\\
t_2^{-1} & 1-t_2^{-2} &\matindex{(8,6)}
\\
\end{block}
\end{blockarray}
\qquad
 -t_2^{-2}
\begin{blockarray}{ccc}
\matindex{(7,8)}
&
\matindex{(8,7)}\\
\begin{block}{(cc)c}
0 & t_1^{-1}&\matindex{(7,8)}
\\
t_1^{-1} & 1-t_1^{-2} &\matindex{(8,7)}
\\
\end{block}
\end{blockarray}
\end{gather*}

\newpage

\noindent
Six 6-dimensional invariant subspaces:
\begin{gather*}
\begin{blockarray}{ccccccc}
\matindex{(1,4)}
&
\matindex{(1,5)}
&
\matindex{(2,3)}
&
\matindex{(3,2)}
&
\matindex{(5,1)}
&
\matindex{(4,1)}\\
\begin{block}{(cccccc)c}
0&0&0&0&0&t_1t_2&\matindex{(1,4)}
\\
0&0&0&0&t_1t_2&0&\matindex{(1,5)}
\\
0&0&0&-\zeta t_1t_2&\zeta t_2\big(1-t_1^2\big)&t_2\big(t_1^2+1\big)&\matindex{(2,3)}
\\[1mm]
0&0&-\zeta t_1t_2&0&t_1\big(t_2^2+1\big)&\zeta t_1\big(1-t_2^2\big)&\matindex{(3,2)}
\\[1mm]
0&t_1t_2&0&t_1\big(t_2^2-1\big)&t_2^2\big(t_1^2-1\big)&\zeta\big(t_2^2-1\big)&\matindex{(5,1)}
\\[1mm]
t_1t_2&0&t_2\big(t_1^2-1\big)&0&\zeta \big(t_1^2-1\big)&t_1^2(t_2^2-1)&\matindex{(4,1)}
\\
\end{block}
\end{blockarray}
\\[-3mm]
\begin{blockarray}{ccccccc}
\matindex{(1,6)}
&
\matindex{(2,4)}
&
\matindex{(2,5)}
&
\matindex{(5,2)}
&
\matindex{(4,2)}
&
\matindex{(6,1)}\\
\begin{block}{(cccccc)c}
0&0&0&0&0&t_2&\matindex{(1,6)}
\\
0&0&0&0&\zeta t_2&\zeta t_2\big(t_1-t_1^{-1}\big)&\matindex{(2,4)}
\\[1mm]
0&0&0&\zeta t_2&0&t_2\big(t_1+t_1^{-1}\big)&\matindex{(2,5)}
\\[1mm]
0&0&\zeta t_2&-\big(t_2^2+1\big)&\zeta\big(t_2^2-1\big)&\zeta \big(t_1t_2^2+t_1^{-1}\big)
&\matindex{(5,2)}
\\[1mm]
0&\zeta t_2&0&0&0&\big(t_1t_2^2+t_1^{-1}\big)&\matindex{(4,2)}
\\[1mm]
t_2&0&t_2\big(t_1-t_1^{-1}\big)&\zeta \big(t_1-t_1^{-1}\big)&t_1\big(t_2^2-1\big)
&\big(t_2^2+t_1^{-2}\big)\big(t_1^2-1\big)&\matindex{(6,1)}
\\
\end{block}
\end{blockarray}
\\[-3mm]
\begin{blockarray}{ccccccc}
\matindex{(1,7)}
&
\matindex{(3,5)}
&
\matindex{(3,4)}
&
\matindex{(4,3)}
&
\matindex{(5,3)}
&
\matindex{(7,1)}\\
\begin{block}{(cccccc)c}
0&0&0&0&0&t_1&\matindex{(1,7)}
\\
0&0&0&0&\zeta t_1&\zeta t_1\big(t_2-t_2^{-1}\big)&\matindex{(3,5)}
\\[1mm]
0&0&0&\zeta t_1&0&t_1\big(t_2+t_2^{-1}\big)&\matindex{(3,4)}
\\[1mm]
0&0&\zeta t_1&-\big(t_1^2+1\big)&\zeta\big(t_1^2-1\big)&\zeta \big(t_1^2t_2+t_2^{-1}\big)
&\matindex{(4,3)}
\\[1mm]
0&\zeta t_1&0&0&0&\big(t_1^2t_2+t_2^{-1}\big)&\matindex{(5,3)}
\\[1mm]
t_1&0&t_1\big(t_2-t_2^{-1}\big)&\zeta \big(t_2-t_2^{-1}\big)&t_2\big(t_1^2-1\big)
&\big(t_1^2+t_2^{-2}\big)\big(t_2^2-1\big)&\matindex{(7,1)}
\\
\end{block}
\end{blockarray}
\\[-3mm]
\resizebox{0.95\textwidth}{!}{
$\begin{blockarray}{ccccccc}
\matindex{(2,8)}
&
\matindex{(5,6)}
&
\matindex{(4,6)}
&
\matindex{(6,4)}
&
\matindex{(6,5)}
&
\matindex{(8,2)}\\
\begin{block}{(cccccc)c}
0&0&0&0&0&-t_1^{-1}&\matindex{(2,8)}
\\[1mm]
0&0&0&0&-\zeta t_1^{-1}&\zeta t_1^{-1}\big(t_2^{-1}-t_2\big)&\matindex{(5,6)}
\\[1mm]
0&0&0&-\zeta t_1^{-1}&0&0&\matindex{(4,6)}
\\[1mm]
0&0&-\zeta t_1^{-1}&0&\zeta \big(t_1^{-2}-1\big)&\zeta t_2^{-1}\big(1-t_1^{-2}\big)
&\matindex{(6,4)}
\\[1mm]
0&-\zeta t_1^{-1}&0&0&-t_1^{-2}\big(1+t_1^{-2}\big)&t_2^{-1}-t_2&\matindex{(6,5)}
\\[1mm]
-t_1^{-1}&\zeta t_1^{-2}\big(t_2^{-1}-t_2\big)&t_1^{-1}\big(t_2^{-1}-t_2\big)
&\zeta \big(t_1^2t_2+t_2^{-1}\big)&-\big(t_1t_2+t_1^{-1}t_2^{-1}\big)
&\big(1+t_1^{-2}t_2^{-2}\big)\big(1-t_2^2\big)&\matindex{(8,2)}
\\
\end{block}
\end{blockarray}$}
\\[-3mm]
\resizebox{0.95\textwidth}{!}{
$\begin{blockarray}{ccccccc}
\matindex{(3,8)}
&
\matindex{(4,7)}
&
\matindex{(5,7)}
&
\matindex{(7,5)}
&
\matindex{(7,4)}
&
\matindex{(8,3)}\\
\begin{block}{(cccccc)c}
0&0&0&0&0&-t_2^{-1}&\matindex{(3,8)}
\\[1mm]
0&0&0&0&-\zeta t_2^{-1}&\zeta t_2^{-1}\big(t_1^{-1}-t_1\big)&\matindex{(4,7)}
\\[1mm]
0&0&0&-\zeta t_2^{-1}&0&0&\matindex{(5,7)}
\\[1mm]
0&0&-\zeta t_2^{-1}&0&\zeta \big(t_2^{-2}-1\big)&\zeta t_1^{-1}\big(1-t_2^{-2}\big)
&\matindex{(7,5)}
\\[1mm]
0&-\zeta t_2^{-1}&0&0&-t_2^{-2}\big(t_2^2+1\big)&\big(t_1^{-1}-t_1\big)&\matindex{(7,4)}
\\[1mm]
-t_2^{-1}&\zeta t_2^{-2}\big(t_1^{-1}-t_1\big)&t_2^{-1}\big(t_1^{-1}-t_1\big)
&\zeta \big(t_1t_2^2+t_1^{-1}\big)&-\big(t_1t_2+t_1^{-1}t_2^{-1}\big)
&\big(1+t_1^{-2}t_2^{-2}\big)\big(1-t_1^2\big)&\matindex{(8,3)}
\\
\end{block}
\end{blockarray}$}
\\[-3mm]
\resizebox{0.95\textwidth}{!}{
$\begin{blockarray}{ccccccc}
\matindex{(4,8)}
&
\matindex{(5,8)}
&
\matindex{(6,7)}
&
\matindex{(7,6)}
&
\matindex{(8,5)}
&
\matindex{(8,4)}\\
\begin{block}{(ccc@{\,\,}ccc)c}
0&0&0&0&0&t_1^{-1}t_2^{-1}&\matindex{(4,8)}
\\[1mm]
0&0&0&0&t_1^{-1}t_2^{-1}&0&\matindex{(5,8)}
\\[1mm]
0&0&0&-\zeta t_1^{-1}t_2^{-1}&\zeta t_2^{-1}\big(t_1^{-2}-1\big)&0&\matindex{(6,7)}
\\[1mm]
0&0&-\zeta t_1^{-1}t_2^{-1}&0&0&\zeta t_1^{-1}\big(t_2^{-2}-1\big)&\matindex{(7,6)}
\\[1mm]
0&t_1^{-1}t_2^{-1}&-\zeta t_1^{-2}\big(t_2+t_2^{-1}\big)&t_1^{-1}\big(1-t_2^{-2}\big)
&-\zeta t_1^{-1}\big(1+t_2^{-2}\big)&\zeta \big(t_2^{-2}-1\big)&\matindex{(8,5)}
\\[1mm]
t_1^{-1}t_2^{-1}&0&t_2^{-1}\big(1-t_1^{-2}\big)&-\zeta t_2^{-2}\big(t_1+t_1^{-1}\big)
&\zeta \big(t_1^{-2}-1\big)&-\zeta t_2^{-1}\big(1+t_1^{-2}\big)&\matindex{(8,4)}
\\
\end{block}
\end{blockarray}$}
\end{gather*}

\newpage

\noindent
One 10-dimensional invariant subspace:
\[
\resizebox{0.99\textwidth}{!}{
$
\begin{blockarray}{@{}c@{\,}c@{\,}c@{\,}c@{\,}c@{\,}c@{\,}c@{\,}c@{\,}c@{\,}c@{\,}c@{}}
\matindex{(1,8)}
&
\matindex{(2,7)}
&
\matindex{(3,6)}
&
\matindex{(4,5)}
&
\matindex{(4,4)}
&
\matindex{(5,5)}
&
\matindex{(5,4)}
&
\matindex{(6,3)}
&
\matindex{(7,2)}
&
\matindex{(8,1)}\\
\begin{block}{(@{}c@{\,}c@{\,}c@{\,}c@{\,}c@{\,}c@{\,}c@{\,}c@{\,}c@{\,}c@{})@{\,\,\,\,\,}c@{}}
0&0&0&0&0&0&0&0&0&1&\matindex{(1,8)}
\\
0&0&0&0&0&0&0&0&1&t_1-t_1^{-1}&\matindex{(2,7)}
\\
0&0&0&0&0&0&0&1&0&t_2-t_2^{-1}&\matindex{(3,6)}
\\[1mm]
0&0&0&0&0&0&-1&\zeta \big(t_1+t_1^{-1}\big)&0&
\big(t_2-t_2^{-1}\big)
\zeta t_1
&\matindex{(4,5)}
\\[1mm]
0&0&0&0&-1&0&0&t_1^{-1}-t_1&0&
\big(t_1-t_1^{-1}\big)
t_2^{-1}
&\matindex{(4,4)}
\\[1mm]
0&0&0&0&0&-1&0&0&t_2^{-1}-t_2&
\big(t_2-t_2^{-1}\big)t_1^{-1}
&\matindex{(5,5)}
\\[1mm]
0&0&0&-1&0&0&0&0&\zeta(t_2+t_2^{-1})&
\big(t_1-t_1^{-1}\big)
\zeta t_2
&\matindex{(5,4)}
\\[1mm]
0&0&1&\zeta \big(t_1+t_1^{-1}\big)&0&t_1-t_1^{-1}&0&0&t_1t_2+t_1^{-1}t_2^{-1}&
\begin{matrix}
\big(t_1t_2+t_1^{-1}t_2^{-1}\big)\\
\times
\big(t_1-t_1^{-1}\big)
\end{matrix}&\matindex{(6,3)}
\\[1mm]
0&1&0&0&t_2-t_2^{-1}&0&\zeta \big(t_2+t_2^{-1}\big)&t_1t_2+t_1^{-1}t_2^{-1}&0&
 \begin{matrix}
\big(t_1t_2+t_1^{-1}t_2^{-1}\big)\\
\times
\big(t_2-t_2^{-1}\big)
\end{matrix}&\matindex{(7,2)}
\\[1mm]
1&t_1-t_1^{-1}&t_2-t_2^{-1}&
\begin{matrix}
 \big(t_2-t_2^{-1}\big)\\
\times \zeta t_1^{-1}
\end{matrix}&
\begin{matrix}
 \big(t_2-t_2^{-1}\big)\\
\times t_1
\end{matrix}
&
\begin{matrix}
\big(t_1-t_1^{-1}\big)\\
\times t_2
\end{matrix}&
\begin{matrix}
\big(t_1-t_1^{-1}\big)\\
\times \zeta t_2^{-1}
\end{matrix}&
\begin{matrix}
\big(t_1t_2+t_1^{-1}t_2^{-1}\big)\\
\times \big(t_1-t_1^{-1}\big)
\end{matrix}&
\begin{matrix}
\big(t_1t_2+t_1^{-1}t_2^{-1}\big)\\
\times \big(t_2-t_2^{-1}\big)
\end{matrix}&
\begin{matrix}
\big(t_1-t_1^{-1}\big)\big(t_2-t_2^{-1}\big)\\
\times
\bigl(t_1t_2+t_1^{-1}t_2^{-1}\bigr)
\end{matrix}&\matindex{(8,1)}
\\
\end{block}
\end{blockarray}
$
}
\]

\subsection*{Acknowledgements}

The authors wish to thank Rinat Kashaev for his generous sharing of his ideas, as well as the anonymous referees for their valuable remarks. MH was partially supported through the NSF-RTG grant
\#DMS-2135960. BMK was partially supported through the BJNSF grant IS24066. GT was supported by the BJNSF grant IS23005.

\pdfbookmark[1]{References}{ref}
\LastPageEnding

\end{document}